\documentclass[twocolumn]{svjour3}          
\smartqed  
\usepackage{graphicx}
\usepackage{custom}
\usepackage{amsmath}
\usepackage{amssymb,url,color,cuted}
\usepackage[numbers]{natbib}
\usepackage{comment}
\usepackage[printwatermark]{xwatermark}
\usepackage{pdfpages}

\usepackage{array}
\usepackage{wrapfig}
\usepackage{multirow}
\usepackage{tabu}
\usepackage{float}
\usepackage[normalem]{ulem}
\useunder{\uline}{\ul}{}
\usepackage{dcolumn}
\usepackage{xcolor}
\usepackage{booktabs}
\usepackage{longtable}
\usepackage{supertabular}
\usepackage{threeparttablex} 
\usepackage{epstopdf}
\usepackage{subfigure}
\usepackage{threeparttablex}
\newcolumntype{V}[1]{>{\centering\arraybackslash}m{#1}}
\newcolumntype{L}[1]{>{\arraybackslash}m{#1}}
\newcolumntype{N}{@{}m{0pt}@{}}

\def\ny#1{{ \color{black}#1}}

\journalname{Health Care Management Science}
\begin{document}

\title{Inverse Optimization on Hierarchical Networks: An Application to Breast Cancer Clinical Pathways
}

\titlerunning{Inverse Optimization on Hierarchical Networks}        

\author{Timothy C. Y. Chan         \and Katharina Forster  \and Steven Habbous\\ Claire Holloway  \and Luciano Ieraci  \and Yusuf Shalaby \and Nasrin Yousefi }

\authorrunning{Chan et al.} 

\maketitle

\begin{abstract}
Clinical pathways are standardized processes that outline the steps required for managing a specific disease. However, patient pathways often deviate from clinical pathways. Measuring the concordance of patient pathways to clinical pathways is important for health system monitoring and informing quality improvement initiatives. 
In this paper, we develop an inverse optimization-based approach to measuring pathway concordance in breast cancer, a complex disease. We capture this complexity in a hierarchical network that models the patient's journey through the health system. A novel inverse shortest path model is formulated and solved on this hierarchical network to estimate arc costs, which are used to form a concordance metric to measure the distance between patient pathways and shortest paths (i.e., clinical pathways). Using real breast cancer patient data from Ontario, Canada, we demonstrate that our concordance metric has a statistically significant association with survival for all breast cancer patient subgroups. We also use it to quantify the extent of patient pathway discordances across all subgroups, finding that patients undertaking additional clinical activities constitute the primary driver of discordance in the population.

\keywords{inverse optimization \and hierarchical network \and clinical pathway concordance \and breast cancer \and survival analysis}
\newpage
\section*{Highlights}
\begin{itemize}
    \item We develop a data-driven metric for measuring the concordance of patient pathways to the clinical pathways for diseases with complex pathways 
    \item We validate the concordance metric by showing its statistically significant association with survival 
    \item The proposed concordance metric can be used to pinpoint the major drivers of discordance and attribute the population discordance to different sections of  the pathways
    \item Quantifying pathway concordance helps policymakers in monitoring the health system and informing quality improvement initiatives 

\end{itemize}

\end{abstract}
\newpage
\section{Introduction}

Clinical pathways are standardized processes that outline the steps required for managing a specific disease \cite{campbell1998integrated,de2006defining}. 
They are designed with the goal of optimizing patient outcomes such as survival or quality of life. Previous research has shown that following clinical pathways can prolong survival, increase patient satisfaction, reduce wait times, reduce in-hospital complication rates, reduce hospital length of stay, and reduce cost of care \cite{rotter2010clinical,rotter2012effects, panella2003reducing,schmidt2018national,vanounou2007deviation,van2013dynamic}. However, patient-traversed pathways can deviate from the recommended pathways for many reasons, and it is important for a health system to be able to identify and quantify the impact of such deviations. Monitoring variation in patient-traversed pathways against the clinical pathways can help decision-makers pinpoint opportunities for improvement in the health system.  
Thus, having tools to facilitate rigorous quantitative measurement of pathway deviation is an important enabler for improving care delivery \cite{forster2020can}. Since ``clinical pathway'' is a general term that can be used to describe different types of products developed to support clinical decision making and streamlining of care, in this paper we use the term \textit{reference pathway} to describe a series of recommended pathways for a group of patients based on their disease characteristics. To give an idea, if the health system is modeled as a network, where nodes represent activities a patient can undertake (e.g., imaging, surgery, etc.), then a reference pathway would be a path in the network.

There exist multiple methods in the literature for measuring the concordance between patient pathways and clinical pathways. The simplest method is to classify patient pathways as ``discordant'' or ``concordant'' by examining if they meet certain pass/fail criteria \cite{cheng2012integrated,bergin2020concordance,forster2020can}. However, there may be subjectivity in choosing which criteria to use. Plus, binary classification of patient pathways ignores the extent of variations in the pathway deviations -- both small and large deviations may be considered discordant, despite the fact that larger deviations may be worse for the patient \cite{karunakaran2020deviations}. Another metric for measuring the similarity of two pathways is edit distance, which treats pathways as strings and measures distance as the number of operations such as additions and deletions needed to transform one string to the other \cite{van2010measuring,yan2018aligning,williams2014using}. 
However, a patient pathway that has an extra activity (addition to reference pathway) may be associated with better outcomes than one that is missing a concordant activity such as treatment (deletion from reference pathway), so a standard, unweighted edit distance metric may not sufficiently differentiate these two cases. A third approach models pathways as walks on a directed graph and uses inverse optimization to learn arc costs that can be used to form a weighted concordance metric based on the cost difference between the walk and a shortest path \cite{chan2021inverse}. For complex diseases though, a graph that models all possible patient journeys in the health system can be very large. 

In this paper, we broaden the applicability of inverse optimization for clinical pathway concordance measurement to diseases with more complex patient journeys by focusing on \emph{hierarchical networks}. When we refer to a hierarchical network, we mean that the overall connected network is composed of nested subnetworks that are defined over multiple ``levels'' -- each subnetwork has nodes that represent subnetworks in a lower level. 
Diseases with complex patient journeys have many different reference pathways which are also long. In addition, the presence or absence of some activities can depend on previous activities. These characteristics are challenging to model using a single network. 
The hierarchical structure facilitates the modeling of long pathways by dividing them into shorter sections, and allows us to incorporate important contextual information by using appropriately designed levels. This structure also helps to model the dependencies of clinical activities by defining the conditions on multiple subnetworks.
We will define the shortest path problem on this hierarchical network in a recursive fashion. The inverse optimization model is derived from the shortest path problem and aims to identify a cost vector such that a given clinical pathway is a shortest path in the hierarchical network.

Inverse optimization aims to infer the parameters of an optimization problem given a set of observed decisions that are assumed to be optimal or suboptimal solutions \cite{troutt2006behavioral, keshavarz2011imputing, chan2014generalized, bertsimas2015data, aswani2018inverse, esfahani2018data, chan2018inverse, babier2021ensemble}. Our problem is an inverse network flow problem, which has been widely studied in the literature \cite{zhang1996network, yang1997inverse, zhang1998inverse, ahuja2001inverse,ahuja2002combinatorial}. The inverse shortest path problem has been studied \cite{burton1992instance, xu1995inverse, zhang1995column, burton1997inverse} and applied to several different application areas including traffic modelling \cite{burton1992instance}, seismic tomography \cite{heuberger2004inverse,burton1992instance}, and healthcare \cite{chan2021inverse}. Most inverse network flow models have been developed for the case of a single observation (i.e., a reference pathway in our context) \cite{zhang1996network, yang1997inverse, zhang1998inverse, ahuja2001inverse,ahuja2002combinatorial}. Some have considered multiple observations \cite{farago2003inverse,zhao2015learning}. See \cite{chan2021review} for a comprehensive review of inverse optimization. Ours is the first to consider multiple observations in a hierarchical network.

Our motivating application is breast cancer. We apply our methodology to measure the concordance of patient pathways against the recommended clinical pathways outlined in the breast cancer pathway maps developed by Ontario Health (Cancer Care Ontario). Pathway maps are essentially flowcharts that provide a high-level view of the care that a cancer patient should receive. 
Ontario Health is the provincial agency responsible for ensuring Ontarians receive high-quality health care services where and when they need them. Ontario Health (Cancer Care Ontario) is an organization within Ontario Health that works to ``equip health professionals, organizations and policy-makers with the most up-to-date cancer knowledge and tools to prevent cancer and deliver high-quality patient care''\citep{ccoabout}. The Disease Pathway Management (DPM) program within Ontario Health (Cancer Care Ontario) ``engages multidisciplinary teams to develop disease site-specific pathway maps that depict the care a typical cancer patient should receive as they progress through the Ontario cancer care system'' \cite{pathdev}. With the goal of optimizing patient outcomes such as survival, these pathway maps are designed for the population and cover the entire continuum of care. 
Breast cancer pathway maps are complex because the recommended treatments vary based on the patient and disease characteristics. 
Moreover, breast cancer treatment consists of multiple different steps or care episodes and the results from previous episodes often determine which is the appropriate next episode along that care continuum.

Breast cancer pathway maps include optional and mandatory clinical activities. Enumerating the reference pathways by considering all the possible combinations of these activities would lead to thousands of reference pathways. 
By using a hierarchical network structure, we divide the entire network into smaller subnetworks and define reference pathways on each subnetwork. This makes the inverse optimization problem tractable by significantly reducing the number of arcs, nodes, and reference pathways. Moreover, some activities in a pathway may depend on previous activities. For example, radiation therapy is recommended for early stage patients who had a breast conserving surgery or if cancer cells are observed at the edge of the removed tissue by surgery, but is discordant otherwise. A single network cannot incorporate this conditional, so one solution is to create more than one network to model each predetermined combination of these conditionals. However, the number of these networks would grow exponentially with the number of such variations in the pathway map. Plus, if concordance scores are calculated using different networks, the scores of the subpopulations will not be comparable. The hierarchical network structure helps with this problem because we can create a level to incorporate these variations in the care network. 

We emphasize that our aim is to develop a metric that measures the concordance of patient pathways against reference pathways, regardless of whether any discordance is justified for the patient at hand. Appropriateness of any pathway deviation would be determined by clinicians using additional data beyond the data needed for concordance measurement. The main use of the concordance metric is to identify opportunities for health system improvement, rather than specific interventions to improve the outcome of a single patient.

Our main contributions in this paper are:
\begin{enumerate}
    \item  Methodology: We propose an inverse optimization approach for learning costs in a hierarchical network that can be used to measure clinical pathway concordance for diseases that have complex pathway maps. Our approach facilitates modeling these complex pathways by sectioning them into meaningful smaller subpathways, applying the methodology to these subpathways, and then aggregating the information from all subpathways to obtain the final concordance score.
    
    \item Application: We apply the proposed methodology to a large dataset of breast cancer patients, who are divided into subgroups based on their disease characteristics and whose journeys are modeled using separate networks. Our numerical results show that:
    \begin{enumerate}
        \item Our concordance metric provides a meaningful quantitative measure of pathway deviation from the clinical pathways by showing that concordance scores for each subgroup have a statistically significant association with survival, even after adjusting for a wide range of confounding variables.
        \item Our metric can pinpoint where subgroup-specific pathway discordances arise. Moreover, the hierarchical network structure facilitates partitioning the discordance of each patient pathway into sections with desired detail. As such, depending on the application, we can examine the population discordance over different levels of the hierarchical network. For example, the higher levels provide a summary view of the discordance while lower levels give more information about the specific clinical activities associated with discordances. 
        
    \end{enumerate}
    
\end{enumerate}

\section{Inverse Optimization Models and Concordance Metric}
\label{sec:iomodel}

In this section, we describe the hierarchical network structure, the inverse optimization models defined on this hierarchical network, and the resulting concordance metric.

\subsection{A Hierarchical Network}

Consider a network with $L+1$ levels. The highest level, level 0, describes the entire network, which we denote $G_0$. We assume there are artificial start ($s_0$) and end ($e_0$) nodes in this network that describe the start and end, respectively, of all patient journeys. We may write these quantities as $G_{0,1}$, $s_{0,1}$ and $e_{0,1}$ when it is notationally convenient. The nodes and arcs between $s_0$ and $e_0$ are partitioned into $T_{1}$ subnetworks $G_{1,1}, \ldots , G_{1,T_{1}}$, which comprise level 1. The subnetworks $G_{1,t}$ have start ($s_{1,t}$) and end ($e_{1,t}$) nodes, and are ordered in the sense that the only arcs entering $s_{1,t}$ come from $s_0$ or $e_{1,u}, u < t$, and the only arcs leaving $e_{1,t}$ go to $s_{1,u}, u > t$, or $e_0$. Level 2 comprises $T_2$ subnetworks $G_{2,1}, \ldots, G_{2,T_2}$. They are partitioned based on their parent subnetwork from level 1, and then ordered within their partition similar to the subnetworks from level 1. That is, $G_{1,1}$ will be composed of subnetworks $G_{2,1}, \ldots , G_{2,T_2^1}$ where $T_2^1 < T_2$, with analogous arc relationships between the start and end nodes of each subnetwork $G_{2,t}, t = 1, \ldots, T_2^1$. Then $G_{1,2}$ will be composed of $T_2^{2}$ subnetworks \\$G_{2,T_2^1+1}, \ldots, G_{2,T_2^1+T_2^2}$ that are connected in a similar fashion, and so on. Finally, each subnetwork $G_{L-1,t}$ at level $L-1$ consists of nodes $1^t_s$, $1^t_e$, $2^t_s$, $2^t_e$, \ldots, $(n_t)^{t}_s$, $(n_t)^{t}_e$ with all arcs of the form $(s_{L-1,t},i^t_s)$, $(i^t_e, e_{L-1,t})$, $(i^t_s,i^t_e)$, $(i^t_e,i^t_s)$, and $(i^t_e,j^t_s)$ for $j > i$. The nodes in $G_{L-1,t}$ represent $n_t$ activities that a patient can undertake in this part of the pathway map. The creation of start ($i^t_s$) and end ($i^t_e$) nodes for each activity $i$ allows us to model the cost of undertaking activity $i^t$ through the cost of the corresponding arc $(i^t_s,i^t_e)$ between these nodes. 

\begin{figure*}
\includegraphics[width=\textwidth]{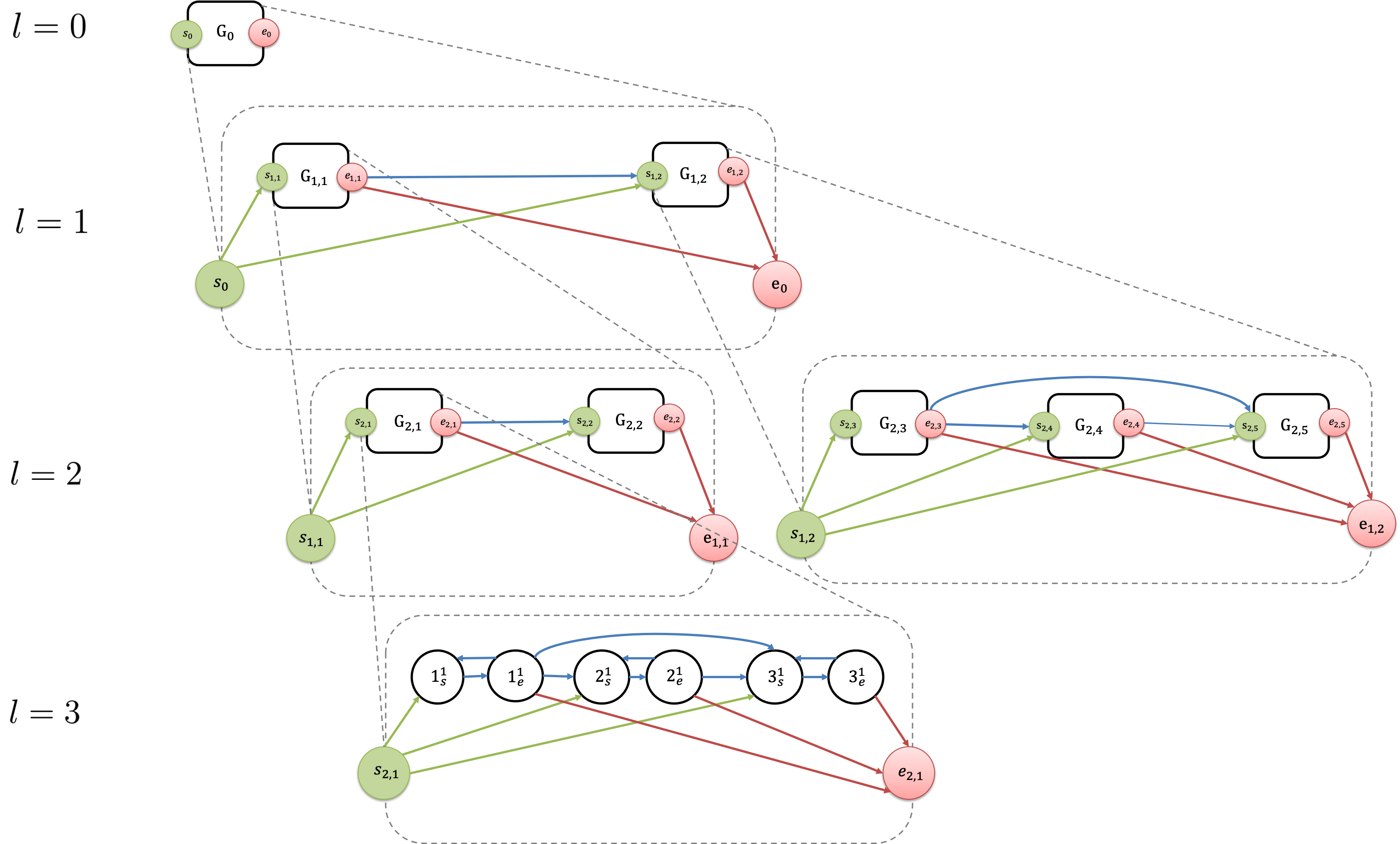}
\caption{An example hierarchical network with $L=3$}
\label{fig:hierarchical}
\end{figure*}

Figure \ref{fig:hierarchical} provides a concrete illustration of a four-level network ($L=3$). The overall network is composed of two subnetworks, $G_{1,1}$ and $G_{1,2}$ in level 1 (i.e., $T_1$ = 2). The first subnetwork $G_{1,1}$ can be further decomposed into $G_{2,1}$ and $G_{2,2}$ (i.e., $T_2^1 = 2$), while the second subnetwork in level 1, $G_{1,2}$, is decomposed into three subnetworks $G_{2,3}$, $G_{2,4}$, and $G_{2,5}$ (i.e., $T_2^2 = 3$). In level 2, $G_{2,1}$ is composed of three activities (i.e., $n_1=3$), represented by $1^1_s$, $1^1_e$, $2^1_s$, $2^1_e$, $3^1_s$, $3^1_e$ as well as the artificial start ($s_{21}$) and end ($e_{21}$) nodes. The other subnetworks in level 2 would have a similar activity-level structure.

The interpretation of our hierarchical network in the concordance measurement application is as follows. Consider a cancer patient with a complex care journey that involves diagnostic activities, pre-surgery activities, surgery activities, and post-surgery activities. The overall network $G_0$ represents the entire healthcare system and a patient journey is a walk through this network. That journey can be divided into the diagnostic, pre-surgery, surgery, and post-surgery phases, which correspond to separate subnetworks in level 1. The diagnostic subnetwork may be further divided into diagnostic imaging and pre-surgery imaging, which are subnetworks in level 2. This process continues until we reach individual diagnostic activities such as ultrasounds and X-rays, which occur in level $L$.  

Finally, note that we design the subnetwoks sequentially meaning that arcs indicate the presence of a subnetwork. For example, the arc between subnetwork $G_{2,3}$ and $G_{2,5}$ in Figure \ref{fig:hierarchical} indicates that $G_{2,3}$ and $G_{2,5}$ were traversed and that $G_{2,4}$ was not traversed. Similarly, we use sequential structure for the activities in level $L$ so that the activity nodes are ordered based on a predefined sequence, and any pathway that is mapped to the graph must be sorted in that sequence.

Before proceeding further, we introduce additional notation for convenience. Let $J_{l,t}$ denote the index set of subnetworks of $G_{l,t}$. 
Let $I_l$ denote the index set of networks in level $l$. That is, $I_l = \cup_t J_{l-1,t}$. Let $\mathcal{A}_{l,t}$ denote the set of arcs in $G_{l,t}$. Meaning that $\mathcal{A}_{l,t}$ includes arcs that connect the subnetworks of $G_{l,t}$. In the next section, we will define reference pathways on subnetwork $G_{l,t}$ that have nonzero values on arcs $\mathcal{A}_{l,t}$.

\subsection{Forward Optimization Model and Hierarchical Reference Pathways}

In this subsection, we formulate the shortest path problem on our hierarchical network. This model serves as the foundation for the inverse optimization model that will be defined in the next subsection. 

Suppose the overall network $G_0$ has $m$ nodes and $n$ arcs. Let $\bc$ and $\bx$ denote the cost vector and flow vector, respectively, on all arcs in the network. A reference or patient pathway is represented as a path or a walk on the network, respectively. Each pathway is represented by a flow vector $\bx$, where $x_{ij}$ is equal to the number of times arc $(i,j)$ is traversed. Therefore, $\bx$ will be a 0-1 vector for a reference pathway but can be greater than 1 for a patient pathway if arc $(i,j)$ is traversed more than once. Let $\bA\in \R^{(m-1)\times n}$ be the node-arc incidence matrix used to describe flow balance constraints on $m-1$ nodes; recall that the last row of the matrix is redundant since flow balance is automatically satisfied at that node if it is satisfied at the rest.

If we formulate the shortest path problem on the entire network $G_0$, disregarding the hierarchical nature of the network, we would have a single linear program $\mathbf{FO}(\bc): \underset{\bx}{\text{min}} \{ \bc'\bx :  \bA \bx = \bb, \bx \ge \bzero\}$ where $\bb \in \R^{m-1}$ is a vector of zeros except with a $+1$ at $s_0$. In this case, the corresponding inverse problem would require the specification of a complete flow vector $\bx$ on all arcs as input. A reference pathway that spans the entire network represented by such a flow vector will be referred to as a \emph{complete} reference pathway. However, for diseases with complex patient journeys, there simply may not be recommended sequences of activities that span the entire network $G_0$, from the start node $s_0$ to the end node $e_0$. In the context of breast cancer specifically, the reference pathways may only specify the recommended guidelines within certain subnetworks of the hierarchical network -- we will refer to these as \emph{hierarchical} reference pathways. Thus, we formulate our forward problem, the shortest path problem, recursively, for each subnetwork $G_{l,t}$. 

Let $z_{l,t}(\bc)$ be the cost of a shortest path in network $G_{l,t}$, given a cost vector $\bc$. For each network $G_{l,t}$, there is a recursive relationship between $z_{l,t}(\bc)$ and $z_{(l+1),u}(\bc)$, $u \in J_{l,t}$. Based on the structure of the network in our hierarchical model, $z_{l,t}(\bc)$ is equal to the cost of the optimal path in network $G_{l,t}$, i.e. the arcs connecting its subnetworks, plus the optimal costs inside the traversed subnetworks. 
Let parameter $\bb_{l,t} \in \R^{m-1}$ be a vector of zeros except with $+1$ at the start node and $-1$ at the end node of network $G_{l,t}$. The matrix $\bD_{l,t} \in \R^{n\times n}$ is a diagonal matrix with 1 in all arc indices in $\mathcal{A}_{l,t}$ and 0 elsewhere. 
We define decision variable $y_{l,u}$ to indicate whether or not subnetwork $G_{l,u}$ is traversed. The recursive forward model for subnetwork $G_{l,t}$ is: 

\begin{equation} \label{model:forwardrec}
\begin{alignedat}{1}
\mathbf{FO}_{l,t}(\bc): \quad \underset{\bx,\by}{\text{minimize}} & \quad \bc'\bx + \sum_{u \in J_{l,t}}{y}_{(l+1),u}z_{(l+1),u}(\bc)\\
\textrm{subject to} & \quad \bA \bx = \bb_{l,t}-\sum_{u \in J_{l,t}} y_{(l+1),u}\bb_{(l+1),u},\\
& \quad (\bI-\bD_{l,t})\bx=0,\\
& \quad \bx \ge \bzero,\\
& \quad y_{(l+1),u} \in \{0,1\}, \quad u \in J_{l,t}.
\end{alignedat}
\end{equation}

The objective function minimizes the total cost of a path through network $G_{l,t}$ by minimizing the cost of the arcs between its subnetworks (first term) plus the cost within those subnetworks (second term). The first constraint enforces flow balance in this hierarchical setting by ensuring that there is a flow surplus of 1 at the start node of $G_{l,t}$ and at the end nodes of $G_{(l+1),u}$ if $y_{(l+1),u}=1$, and a flow deficit of $-1$ at the end node of $G_{l,t}$ and at the start nodes of $G_{(l+1),u}$ if $y_{(l+1),u}=1$. The second constraint ensures zero flow outside of network $G_{l,t}$ (or any of its subnetworks). Note that networks $G_{L-1,t}$ do not have any subnetworks, and the problem reduces to a shortest path problem with one additional constraint (second constraint), which restricts the flow to the arcs between the activities of this subnetwork. 
We can solve the recursive shortest path problems by starting from networks in level $L-1$ and working backwards until we reach level 0. 

To illustrate how the recursive forward problem works, consider network $G_{1,1}$. For a given cost vector $\bc$, a shortest path problem on $G_{1,1}$ perceives the two subnetworks $G_{2,1}$ and $G_{2,2}$ as nodes of $G_{1,1}$. However, these subnetworks are special types of nodes as they incur some cost due to the network structure inside them. Therefore, $\mathbf{FO}_{1,1}(\bc)$ should minimize the total cost of a path in $G_{1,1}$ plus the cost incurred by $G_{2,1}$ and $G_{2,2}$. Thus, the objective function becomes $\bc'\bx + {y}_{2,1}z_{2,1}(\bc) + {y}_{2,2}z_{2,2}(\bc)$ where the model decides whether a shortest path should go through $G_{2,1}$ (${y}_{2,1}=1$, ${y}_{2,2}=0$) only, $G_{2,2}$ (${y}_{2,1}=0$, ${y}_{2,2}=1$) only, or both (${y}_{2,1} = {y}_{2,2}=1$). Now, the flow balance constraints in \eqref{model:forwardrec} enforce that a shortest path starts from $s_{1,1}$ and ends at $e_{1,1}$. In addition, they ensure that flow balance holds for the start and end nodes of the subnetworks that are traversed. For example, suppose that a shortest path goes through subnetwork $G_{2,1}$ only. A unit flow that comes from $s_{1,1}$ must exist from $s_{2,1}$, enter from $e_{2,1}$, and exit from $e_{1,1}$. There is no incoming or outgoing flow at nodes $s_{2,2}$ and $e_{2,2}$. Finally, although $\bx$ denotes the flow vector for all the arcs in the overall network, the shortest path problem for $G_{1,1}$ only focuses on the arcs that connect the subnetworks, i.e.,  $A_{1,1}=\{(s_{1,1},s_{2,1}),(s_{1,1},s_{2,2}),(e_{2,1},s_{2,2}),(e_{2,1},e_{1,1}),(e_{2,2},e_{1,1})\}$. The third constraint in \eqref{model:forwardrec} sets the flow on all the other arcs to 0.

Note that with the recursive forward model, we transform a single large continuous problem into a number of small discrete problems as we set the flow on arcs outside of the corresponding subnetwork to 0. This transformation does not add to the complexity of the inverse problem because the inverse problem aims to impute the cost vector, $\bc$, given a flow vector, and once the flow vector is specified, the discrete variables $\by$ are fixed. And then the added benefit of this approach is that the special structure of the recursive formulation affords the use of hierarchical reference pathways as its inputs. 

\subsection{An Inverse Optimization Model Using Hierarchical Reference Pathways}
\label{sec:hiermodel}

The goal of the inverse optimization model is to identify a cost vector $\bc \in \R^n$ that minimizes the suboptimality (with respect to being shortest paths) of a set of reference pathways. This approach generalizes the classical inverse shortest path problem since if it is possible for all reference pathways to simultaneously be shortest paths for some cost vector, the inverse optimization model will identify such a cost vector.  Next, we present two inverse optimization models that, taken together, form our inverse optimization framework for estimating a cost vector. The first model uses reference pathways only to find a cost vector. The second formulation uses patient-traversed pathways to refine the optimal cost vector returned by the first model, in order to improve the association between pathway concordance and survival.

Let $\hat{\mathcal{X}}_f=\{\hat\bx^{1}_{f}, \ldots, \hat\bx^{R}_{f}\}$ denote a set of $R$ complete reference pathways. Assume these flow vectors specify the flow on each arc in the network. Then, similar to the forward problem, we may formulate the inverse problem without considering the hierarchical network structure. We proceed to do so first, to build the foundation for our inverse model on the hierarchical network, which follows. To find a cost vector that makes a given set of complete reference pathways optimal, we need to write down the optimality conditions associated with the shortest path problem, which we can do easily using linear programming duality. Since the complete reference pathways are assumed to be feasible for the forward problem (paths), all we need is to ensure the cost vector satisfies dual feasibility and strong duality. Instead of enforcing strong duality exactly though, we introduce a duality gap, which is minimized, ensuring that our model is feasible even when the complete reference pathways cannot simultaneously be shortest paths. The following formulation is the inverse optimization model that minimizes the sum of squared absolute duality gaps induced by the cost vector $\bc$ and the set of complete reference pathways on the entire network:
\begin{equation} \label{model:hieriofull}
\begin{alignedat}{1}
\underset{\bc, \bp,\boldsymbol{\epsilon}_f}{\text{minimize}} & \quad \sum_{q=1}^{R}(\epsilon^{q}_{f})^2\\
\textrm{subject to} & \quad \bA' \bp \le \bc,\\
& \quad \bc'\hat\bx^{q}_{f}=\bb'\bp+\epsilon^{q}_{f}, \quad   q=1, \ldots, R, \\
& \quad \|\bc\|_{\infty}= 1,\\
& \quad \bA \bc= \bzero,
\end{alignedat}
\end{equation}

\ny{The objective function minimizes the sum of squared duality gaps rather than absolute ones to avoid large duality gaps}. The vector $\bp \in \R^{m-1}$ is the dual vector associated with flow balance constraints over the entire network. The first constraint represents dual feasibility for the entire network and the second constraint defines the duality gap $\epsilon^{q}_{f}$ for each complete reference pathway $\hat\bx^{q}_{f}$. We add the third and fourth constraints to ensure the inverse model returns a meaningful cost vector. The third constraint ensures the cost vector is not zero. 
The fourth constraint is needed to ensure that $\bc$ lies in the lower dimensional space defined by the span of the columns of $\bA$ \ny{in the network flow problem}, and is not orthogonal to the entire feasible region. An orthogonal cost vector $\bc$ makes all the feasible solutions \ny{of the network flow problem} optimal, i.e. every reference and patient pathway would be optimal. We want to avoid this case and find a cost vector that can differentiate between different pathways and generate meaningful discordance scores. 

As a reminder, this formulation is applicable if the reference pathways are complete, i.e., they specify the flow for every arc in the network. However, in general, reference pathways are only specified for a subset of the network. In principle, we can form complete reference pathways by combining the hierarchical reference pathways from different subnetworks, including all combinations across all levels, but this would result in an explosion in the number of complete reference pathways. Therefore, we incorporate the hierarchical network structure into the inverse problem and develop a new inverse optimization formulation tailored for a hierarchical network with hierarchical reference pathways. 

Let $\hat{\mathcal{X}}_{l,t}=\{\hat\bx^{1}_{l,t}, \ldots, \hat\bx^{R_{l,t}}_{l,t}\}$ denote the set of $R_{l,t}$ hierarchical reference pathways for network $G_{l,t}$, each assumed to be a feasible solution to $\mathbf{FO}_{l,t}(\bc)$, i.e., a path. Each $\hat{\bx}^{q}_{l,t} \in \R^n$ is a disjoint flow vector which has nonzero flow only for arcs in $\mathcal{A}_{l,t}$ and carries information on the sets of subnetworks that are traversed. Let $K^{q}_{l,t} \in J_{l,t}$ denote the set of subnetworks that are traversed by $\hat{\bx}^{q}_{l,t}$. Finally, let $\hat{\mathcal{X}}= \cup_{l=0}^{L-1} \cup_{t=1}^{|I_l|}  \hat{\mathcal{X}}_{l,t}$ be the set of hierarchical reference pathways over all subnetworks.  
The hierarchical inverse optimization formulation is
\begin{strip}
\begin{equation} 
\label{model:hierioref}
\begin{alignedat}{1}
\mathbf{IO}^{\textrm{ref}}(\hat{\mathcal{X}}): \quad \underset{\bc, \bp,\boldsymbol{\epsilon}}{\text{minimize}} & \quad \sum_{l=0}^{L-1}\sum_{t=1}^{|I_l|}\sum_{q=1}^{R_{l,t}}(\epsilon^{q}_{l,t})^2\\
\textrm{subject to} & \quad \bA' \bp \le \bc,\\
& \quad \bc'\hat\bx^{q}_{l,t}+\sum_{u\in K^{q}_{l,t}}\bb_{(l+1),u}^\prime\bp=\bb'_{l,t}\bp+\epsilon^{q}_{l,t}, \quad   q=1, \ldots, R_{l,t}, t \in I_l, l\in\{0, \ldots, L-1\}, \\
& \quad \|\bc\|_{\infty}= 1,\\
& \quad \bA \bc= \bzero.
\end{alignedat}
\end{equation}
\end{strip}
The first, third, and fourth constraints are defined in the same way as the ones in \eqref{model:hieriofull}. We embed the recursive structure of \eqref{model:forwardrec} into the strong duality constraint from \eqref{model:hieriofull}, resulting in separate strong duality constraints for each hierarchical reference pathway over all subnetworks. Thus, the second constraint in \eqref{model:hierioref} defines the duality gap $\epsilon^{q}_{l,t}$ for each hierarchical reference pathway $\hat\bx^{q}_{l,t}$ using the relation between the objective value of the forward problem on $G_{l,t}$ and its subnetworks. 
The first term $\bc'\hat\bx^{q}_{l,t}$ is the cost of flow on the arcs between each traversed subnetwork. As described in formulating the forward problem, when we focus on network $G_{l,t}$, its subnetworks are seen as special nodes that come with their own cost. The hierarchical reference pathways defined on $G_{l,t}$ do not include any information about the hierarchical reference pathways within each subnetwork, so we explicitly add the optimal cost within these subnetworks as $\sum_{u\in K^{r}_{l,t}}\bb_{(l+1),u}^\prime\bp$. 
Finally, $\bb'_{l,t}\bp$ is the shortest path cost through network $G_{l,t}$. 
Note that each duality gap variable is nonnegative since the hierarchical reference pathways are feasible for the forward problem in each network.

In the rest of this paper, we use the inverse optimization model \eqref{model:hierioref}. The next subsection explains how we further refine the cost vector.

\subsection{Refining the Solution Using Patient Data}
\label{sec:hiermodel_patient}

Let $\hat{\mathcal{W}}_{l,t}=\{\hat\bw^{1}_{l,t}, \ldots, \hat\bw^{S_{l,t}}_{l,t}\}$ denote a dataset of $S_{l,t}$ patient pathways with positive clinical outcomes (i.e., survived) and $\hat{\mathcal{V}}_{l,t}=\{\hat\bv^{1}_{l,t}, \ldots, \hat\bv^{D_{l,t}}_{l,t}\}$ denote a dataset of $D_{l,t}$ patient pathways with negative clinical outcomes (i.e., died) for network $G_{l,t}$. All of these pathways are assumed to be feasible flow vectors for the shortest path problem on their corresponding networks. Let $\hat{\mathcal{W}}= \cup_{l=0}^{L-1} \cup_{t=1}^{|I_l|} \hat{\mathcal{W}}_{l,t}$ and $\hat{\mathcal{V}}= \cup_{l=0}^{L-1} \cup_{t=1}^{|I_l|} \hat{\mathcal{V}}_{l,t}$ be the set of patient pathways with positive and negative clinical outcomes over all the subnetworks in the model.
Let $M^{q}_{l,t} \in J_{l,t}$ and $N^{q}_{l,t} \in J_{l,t}$ be the index sets of networks traversed by $\hat\bw^{q}_{l,t}$ and $\hat\bv^{q}_{l,t}$, respectively.
If formulation~\eqref{model:hierioref} has multiple optimal solutions, we can use patient data to determine a cost vector that, in addition to maximizing fit with the hierarchical reference pathways, provides separation between $\hat{\mathcal{W}}$ and $\hat{\mathcal{V}}$. Once formulation~\eqref{model:hierioref} is solved and the set of optimal duality gaps is generated, we use it as input into the following inverse optimization model: 
\newpage 
\begin{strip}
\begin{equation} 
\label{model:hierpatientref}
\begin{alignedat}{1}
\mathbf{IO}^{\textrm{pat}}(\hat{\mathcal{W}},\hat{\mathcal{V}},\boldsymbol{\epsilon}^*): \quad \underset{\bc, \bp, \boldsymbol{\gamma}, \boldsymbol{\delta}}{\text{minimize}} & \quad \sum_{l=0}^{L-1} \alpha_{l}\sum_{t=1}^{|I_l|} \sum_{q=1}^{S_{l,t}} \gamma^{q}_{l,t}-\sum_{l=0}^{L-1} \beta_{l} \sum_{t=1}^{|I_l|} \sum_{q=1}^{D_{lt}} \delta^{q}_{l,t}\\
\textrm{subject to} & \quad \bA' \bp \le \bc,\\
& \quad \bc'\hat\bx^{q}_{l,t}+\sum_{u\in K^{q}_{l,t}}\bb_{(l+1),u}^\prime\bp=\bb_{l,t}'\bp+\epsilon^{q*}_{l,t}, \quad q=1, \ldots, R_{l,t}, t \in I_l, l\in\{0, \ldots, L-1\},   \\
& \quad \bc'\hat\bw^{q}_{l,t}+\sum_{u\in M^{q}_{l,t}}\bb_{(l+1),u}^\prime\bp=\bb_{l,t}'\bp+\gamma^{q}_{l,t}, 
\quad q=1, \ldots, S_{l,t}, t \in I_l,  l\in\{0, \ldots, L-1\},  \\
& \quad \bc'\hat\bv^{q}_{l,t}+\sum_{u\in N^{q}_{l,t}}\bb_{(l+1),u}^\prime\bp=\bb_{l,t}'\bp+\delta^{q}_{l,t}, \quad q=1, \ldots, D_{l,t}, t \in I_l,  l\in\{0, \ldots, L-1\},  \\
& \quad \|\bc\|_{\infty}= 1,\\
& \quad \bA \bc= \bzero.
\end{alignedat}
\end{equation}
\end{strip}

The second constraint forces the cost vector to achieve the optimal duality gap found in model~\eqref{model:hierioref}. In other words, all feasible solutions of~\eqref{model:hierpatientref} are optimal solutions to~\eqref{model:hierioref}. The third and fourth constraints define the duality gaps with respect to the pathways in $\hat{\mathcal{W}}_{l,t}$ and $\hat{\mathcal{V}}_{l,t}$, respectively. Again, the duality gaps are nonnegative since all patient pathways are feasible for the forward problem. The remaining constraints are as defined in~\eqref{model:hierioref}. Unlike model~\eqref{model:hierioref}, the duality gaps in the objective function are linear, as quadratic terms would result in a non-convex problem.
Note that the number of patient pathways for different levels may vary. Moreover, the number of patients with positive and negative outcomes can be disproportionate. To compensate for the imbalance in the data for different levels and patient outcomes, we use the parameters $\alpha_{l}$ and $\beta_{l}$ in the objective function. They are defined as follows. Let $S_l$ be the total number of patient pathways with positive outcomes in level $l$, i.e., $S_l = \sum_{t\in I_l} S_{l,t}$. Similarly, let $D_l = \sum_{t\in I_l} D_{l,t}$. We define $B = \underset{l \in \{0,\ldots,L-1\}}{\text{min}} \{ S_l + D_l\}$ and set $\alpha_l = B/S_l$ and $\beta_l = B/D_l$. 

Note that model~\eqref{model:hierpatientref} is essential in finding a meaningful cost vector in our inverse optimization framework. The reference pathways used in model~\eqref{model:hierpatientref} can be a small subset of the possible pathways and may not cover all the arcs in the subnetworks. Adding the patient data in model~\eqref{model:hierpatientref} helps to find an optimal cost vector among all the optimal cost vectors from model~\eqref{model:hierioref} that will result in a more meaningful measure of concordance. \ny{In our numerical results, we show that without model~\eqref{model:hierpatientref}, the association between the concordance metric and survival may become statistically insignificant.}

\subsection{Concordance Metric}
\label{sec:metric}

In this section, we explain how we measure the concordance of a given patient pathway against the hierarchical reference pathways in a hierarchical network. First, we need to map each patient pathway to the subnetworks. The mapping is a top-down procedure meaning that we start from level 0 and recursively reach level $L-1$. Then, we explain how we use the recursive structure to calculate the discordance of a patient pathway.  Finally, we use the discordance to construct a concordance metric.

Let $\hat\bx_{l,t}$ be a walk on network $G_{l,t}$ corresponding to patient pathway $\hat\bx$ and let $\hat K_{l,t}\in J_{l,t}$ be the index set of subnetworks that are traversed by $\hat\bx_{l,t}$. 
We define the myopic discordance of $\hat\bx$ through network $G_{l,t}$ as 
\begin{equation}
\label{eq:epsilonm}
    \bar\epsilon_{l,t} (\hat\bx)= \bc^*{'}\hat\bx_{l,t}+\sum_{u\in \hat K_{l,t}}\bb_{(l+1),u}^{'}\bp^*-\bb_{l,t}^{'}\bp^*,
\end{equation}
where $\bc^*$ and $\bp^*$ are the optimal cost vector and dual vector generated from the inverse optimization model. The myopic discordance is the optimality gap for the patient pathway on network $G_{l,t}$ alone, which is adapted from the strong duality constraints in the inverse problem. 

In contrast to the myopic discordance, we define the aggregated discordance as the sum of the optimality gap of the patient pathway on network $G_{l,t}$ plus its subnetworks. 
We first define the aggregated discordance for level $L-1$. Recall that level $L-1$ does not contain any subnetworks and consists of the clinical activities. Thus, the second term in the above equation is not relevant and the aggregated discordance of $\hat\bx$ through network $G_{L-1,t}$ is
\begin{equation}
\label{eq:epsilonl1}
    \epsilon_{L-1,t} (\hat\bx) = \bc^*{'}\hat\bx_{L-1,t} -\bb_{L-1,t}^{'}\bp^*.
\end{equation}
The aggregated discordance through networks in all the other levels is calculated using backward recursion
\begin{equation}
\label{eq:epsilona}
    \epsilon_{l,t} (\hat\bx) =
    \bar\epsilon_{l,t} (\hat\bx) + 
    \sum_{u\in \hat K_{l,t}}\epsilon_{(l+1),u} (\hat\bx), \quad l \in \{0, \ldots, L-2\} 
\end{equation}

This equation states that the total discordance of a patient pathway $\hat\bx$ through network $G_{l,t}$ is due to the subnetworks that are traversed (first term) and the discordance within those subnetworks (second term). This split is useful because it allows us to precisely identify the sources of discordance for any pathway in a later analysis (Section \ref{sec:pod}). 

To calculate the total discordance of patient pathway $\hat\bx$ on the entire network $G_0$, we first start from level $L-1$ and use equation \eqref{eq:epsilonl1} to find the optimality gap for networks in level $L-1$ that are traversed by $\hat\bx$. Then, we calculate the optimality gaps for their parent networks until we reach level 0 with backward recursion. The total discordance of patient pathway $\hat\bx$ is equal to $\epsilon_{0,1}(\hat\bx)$, which we write as $\epsilon(\hat\bx)$. 

To illustrate, suppose that we are given an arbitrary patient pathway $\hat\bx^0$ and we map it to the hierarchical network in Figure \ref{fig:hierarchical}.
Suppose that $\hat\bx^0$ traverses only $G_{1,1}$ in level 1 and $G_{2,1}$ and $G_{2,2}$ in level 2. We first calculate $\epsilon_{2,1}(\hat\bx^0)$ and $\epsilon_{2,2}(\hat\bx^0)$. Then, for $G_{1,1}$, we find the aggregate optimality gap, which is the sum of its myopic optimality gap and the aggregated optimality gaps of $G_{2,1}$ and $G_{2,2}$, i.e., $\epsilon_{1,1}(\hat\bx^0)= \bar\epsilon_{1,1}(\hat\bx^0)+\epsilon_{2,1}(\hat\bx^0)+\epsilon_{2,2}(\hat\bx^0)$. Finally, the total optimality gap of $\hat\bx^0$ is $\epsilon(\hat\bx^0) := \epsilon_{0,1}(\hat\bx^0)= \bar\epsilon_{0,1}(\hat\bx^0)+\epsilon_{1,1}(\hat\bx^0)$.

Finally, the concordance metric $\omega(\hat\bx)$ is defined using the discordance as 
\begin{equation}
\label{eq:omega}
    \omega(\hat\bx) = 1 - \frac{\epsilon(\hat\bx)}{\Lambda(\hat\bx)}, 
\end{equation}
where $\Lambda(\hat\bx)$ is the total discordance of a walk that has the largest discordance among all the walks with steps up to the number of steps in $\hat\bx$.  Calculating $\Lambda(\hat\bx)$ can be done by subtracting the cost of a shortest path from the cost of this walk. The corresponding walk and the walk cost are calculated using dynamic programming \cite{chan2021inverse}.  Normalizing the concordance metric in this way ensures that the concordance metric is between 0 and 1. 

\section{Application to Breast Cancer}
\label{sec:application}

In this section, we demonstrate the application of our framework to measuring clinical pathway concordance in breast cancer.

\subsection{Breast Cancer and Subgroups}
Breast cancer is a malignant tumor that originates in the breast. According to recent statistics, 
it is the most commonly diagnosed cancer in the world and accounts for 12\% of new cancer cases \cite{whobreastdata}. 
It is the second leading cause of cancer death in women \cite{statbreast}. Breast cancer can occur in both men and women, but more than 99\% of breast cancer patients are women \cite{statbreast2}. 

Our work is based on the pathway maps used in Ontario, Canada. The Disease Pathway Management (DPM) program at Ontario Health (Cancer Care Ontario) works with multidisciplinary groups of clinicians with disease site expertise to develop pathway maps that outline evidence-based best practices for managing various cancer types through all phases of cancer care, from screening through follow-up care \citep{pathdev}. In this paper, we evaluate concordance with the screening, diagnosis, and treatment portions of the breast cancer care continuum (see pathway maps in \cite{breastpathways}; we use version 2015.10).
A pathway map shows all the possible ways a patient should traverse the care network. Since the clinical pathways differ based on disease characteristics and information that becomes available along the way, the pathway maps are large and complex, containing multiple decision and branch points along the trajectory. Key factors in treatment decision making are cancer stage and biomarker status. Cancer staging includes the T stage, which describes the size of the tumor, the N stage, which indicates the spread of the cancer to the surrounding lymph nodes, and the M stage, which indicates presence/absence of metastases. Biomarker status indicates the expression of the estrogen receptor (classified as ER-positive) and overexpression of the human epidermal growth factor receptor 2 (HER2) (classified as HER2-positive). 
We group patients into 10 subgroups based on TNM stage and biomarker status. Table \ref{tab:subgroups} shows the definitions of these 10 subgroups.

\begin{table*}
\caption{Patient subgroups}
\label{tab:subgroups}
\begin{tabular}{cccccc}
\toprule
Subgroup                                                      & T stage                                                             & N stage                                          & M stage                                       & ER                                            & HER2                                                \\
\toprule
1        & Tis                                                             & 0                                                                                                               & 0                                             &   n/a                                                 &    n/a                                                \\ \hline
2        & 1,2,3                                                               & 0                                                & 0                                             & Pos.                                                & Neg.                                                \\ \hline
3            & 0,1,2                                                               & 1                                                & 0                                             & Pos.                                                & Neg.                                                \\ \hline
4                                                           & \begin{tabular}[c]{@{}c@{}}0,1,2\\ 3\end{tabular}     & \begin{tabular}[c]{@{}c@{}}0,1\\ 0\end{tabular}  & 0 & Neg. &  Neg. \\ \hline
5                                                       & 1a,1b                                                               & 0                                                & 0                                             &  any                                                & Pos.                                                \\ \hline
6                                                           & \begin{tabular}[c]{@{}c@{}}0,1,2\\ 3\end{tabular}     & \begin{tabular}[c]{@{}c@{}}0,1 \\ 0\end{tabular} & 0 &   any  & Pos. \\ \hline
7                                                           & \begin{tabular}[c]{@{}c@{}}4\\ any \\3\end{tabular}     & \begin{tabular}[c]{@{}c@{}} any\\ 2,3\\ 1\end{tabular} & 0 &   any &  Pos. \\ \hline
8                                                                 & \begin{tabular}[c]{@{}c@{}}4\\ any \\3\end{tabular}     & \begin{tabular}[c]{@{}c@{}} any \\ 2,3\\ 1\end{tabular} & 0 & Neg.   &  Neg. \\ \hline
9                                                                 & \begin{tabular}[c]{@{}c@{}}4\\ any\\3\end{tabular}     & \begin{tabular}[c]{@{}c@{}} any\\ 2,3\\ 1\end{tabular} & 0 & Pos.   &  Neg. \\ \hline
10                                                                  & any                                                                 & any                                                 &         1                                      &   any                                               &   any  \\   
\bottomrule
\end{tabular}
    \begin{tablenotes}
      \scriptsize
      \item   \textbf{Tis}: carcinoma in situ, \textbf{T0}: no evidence of cancer in the breast,
      \item  \textbf{T1a}: tumor size is between 1 and 5 mm, \textbf{T1b}: tumor size is between 5 and 10 mm,
      \item \textbf{T2}: tumor size is between 20 and 50 mm, \textbf{T3}: tumor is larger than 50 mm,
      \item  \textbf{T4}: tumor has grown into the chest wall or skin or both.
      \item  \textbf{N0}: no cancer / areas of cancer smaller than 0.2 mm in the lymph nodes, 
      \item  \textbf{N1}: cancer spread to 1-3 axillary lymph nodes
      \item  \textbf{N2}: cancer spread to 4-9 axillary lymph nodes, \textbf{N3}: cancer spread to 10+ axillary lymph nodes.
      \item  \textbf{M0}: no distant metastases, \textbf{M1}: metastasis to another part of the body. 
    \end{tablenotes}
\end{table*}

Each subgroup is part of the same pathway map, but has a different set of reference pathways, since the treatment recommendations depend on the characteristics that define the subgroup. In our analysis, we exclude subgroup 10 since patients in this subgroup are palliative, and care often includes activities not reflected in the pathway maps. A patient might change subgroups in some cases, but our data includes the most recent subgroup, and we do not have information about the exact transition point. This is a limitation of our data set. However, the subgroup change often occurs around the surgery time, and since the pathways are similar before the surgery, the concordance score estimation error should be minimal.

\subsection{Data}
We obtained data of patients in Ontario who had a new breast cancer diagnosis between 2010 and 2016. Their clinical activities were recorded from as early as 2006 until as late as January 2020 or death. We excluded patients who had a missing or invalid health insurance number. Patients who received a clinical diagnosis only, had no treatment or other healthcare activities, or were diagnosed by autopsy or at the time of death were also excluded. We omitted patients with previous cancer (except for basal and squamous cell skin cancers) any time before diagnosis since any prior cancer may alter a patient's diagnosis and treatment pathway. We are unable to definitively identify patients with cancer recurrence due to data limitations. Therefore, we chose to exclude patients with a second cancer (of any type) occurring less than two years after the breast cancer as the diagnosis of a new cancer can affect the treatment follow-ups. We also excluded patients who died within 15 months of diagnosis or did not have all the covariates to be used in survival analysis (See Section \ref{sec:survival}).

Our dataset includes the assigned subgroup for each patient. Since subgroup 1 comprises only 0.0025\% of patients, making it impossible to validate the concordance metric using the survival analysis, we omit this subgroup from our analysis and focus on subgroups 2 to 9. The final cohort contains 47,312 patients, with subgroup 2 being the largest with approximately 48\% of patients. For subsequent analyses, we split the patient data into two sets. The ``training" set contains patients diagnosed in 2010-2012 (approximately 40\% of the data) and the ``testing" set contains patients diagnosed in 2013-2016 (approximately 60\% of the data).  
The training set will be used to learn the arc costs in the network and the testing set will be used for out-of-sample analyses. Table \ref{tab:datasplit} shows the number of data points in each subgroup in the training and testing datasets. 

\begin{table}
\caption{Data split for the model }
\label{tab:datasplit}
\begin{tabular}{cccc}
\toprule
Subgroup & \begin{tabular}[c]{@{}c@{}}Training\\ (n=19,271)\end{tabular} & \begin{tabular}[c]{@{}c@{}}Testing\\ (n=28,041)\end{tabular} & \begin{tabular}[c]{@{}c@{}}Total\\ (n=47,312)\end{tabular} \\
\toprule
2        & 9,339                                                         & 13,718                                                       & 23,057                                                     \\
3        & 3,141                                                        & 4,403                                                        & 7,544                                                      \\
4        & 1,772                                                         & 2,520                                                        & 4,292                                                      \\
5        & 853                                                          & 1,313                                                       & 2,166                                                      \\
6        & 1,306                                                         & 2,198                                                        & 3,504                                                      \\
7        & 689                                                          & 971                                                        & 1,660                                                      \\
8        & 377                                                          & 501                                                         & 878                                                      \\
9        & 1,794                                                         & 2,417                                                       & 4,211    \\
\bottomrule
\end{tabular}
\end{table}

The dataset consists of event logs for each patient describing the clinical activities undertaken (e.g., imaging, treatment, etc.) and their corresponding date stamps. We refine each patient pathway to only include activities from 6 months before diagnosis to 15 months after diagnosis, since patients are expected to complete diagnosis and treatment activities within 15 months of diagnosis. For patients who had multiple surgeries, we only record the first surgery and a follow-up surgery and discard further follow-up surgeries from the pathway in order to simplify the network. These cases were also rare in our dataset.

\subsection{Network Design}
\label{sec:breastnw}
In this section, we describe the high-level pathway requirements for all the included subgroups. We then discuss the detailed pathway requirements of subgroup 2 and explain the step-by-step process of designing the network. We take a similar approach for the other subgroups, but we omit the details for brevity.

Screening and diagnosis, neo-adjuvant therapy, surgery, and adjuvant therapy are the major categories of activities in the breast cancer pathway maps. Screening and diagnosis can include different types of consultations, imaging, biopsies, and staging tests, which are necessary for all subgroups.
Neo-adjuvant therapy may include hormone therapy, chemotherapy, targeted therapy, or a combination of these. Neo-adjuvant therapy is required for subgroups 7, 8, and 9, optional for subgroups 4 and 6, and not required for the others.
Surgery can either be a mastectomy, in which the entire breast is removed, or a breast-conserving surgery (BCS). Pathways in all subgroups must have one breast surgery. Along with the breast surgery, sentinel lymph node biopsy (SLNB) or axillary lymph node dissection (ALND) is required. Subgroups 2 and 5 must have one SLNB, subgroups 4 and 6 must have one SLNB or ALND, and subgroups 3,7,8,9 must have at least one of the two. 
More imaging is often required on the day of the surgery to guide the operation. Adjuvant therapy can include radiation therapy in addition to hormone therapy, chemotherapy, and targeted therapy. Adjuvant therapy is required for all subgroups, but it is optional for subgroups 4 and 6 if they had neo-adjuvant therapy.

We model the breast cancer care network as a hierarchical network with five levels ($L=4$), where level 0 consists of the entire network and level 4 contains the clinical activities. We then build hierarchical reference pathways in each level based on the treatment requirements and map both reference and patient pathways to the resulting network. We use one network structure for all the subgroups, but the reference pathways and the details vary for each subgroup. We explain the network design for subgroup 2 below. 

\subsubsection{Level 1}
The hierarchical network starts with the overall network in level 0. This high-level network consists of five subnetworks in the subsequent level, level 1, that describe different sections of the pathway. These five networks represent:
\begin{enumerate}
    \item Diagnosis: from six months before diagnosis day until the first observed surgery or therapy (non-surgery treatment)
    \item Neo-adjuvant: from the first observed therapy before surgery until the first surgery
    \item Surgery: from the first surgery until the first adjuvant therapy
    \item Adjuvant: from the first therapy after surgery until 15 months after diagnosis
    \item Continual: used for patient pathways in subgroups 4 and 6 to track all targeted therapy and chemotherapy activities throughout the pathway 
\end{enumerate}

We derive the requirements from the breast cancer pathway maps, augmented with clinical domain expertise. To reflect these requirements in the hierarchical network model, we assign `mandatory', `optional', or `discordant' labels to networks in this level. Missing a network with `mandatory' label is considered discordant. In contrast, the presence of a `discordant' network is deemed discordant. Traversing or avoiding an `optional' network does not affect concordance of a given patient pathway. 

For subgroup 2, 
the `Diagnosis', `Surgery', and `Adjuvant' networks are mandatory, while `Neo-adjuvant' and `Continual' are discordant. Hierarchical reference pathways are generated using the combination of mandatory and optional networks, if any. The mandatory ones must be present in all hierarchical reference pathways, while all possible permutations of the optional ones result in a multitude of different hierarchical reference pathways. Thus, \textit{Diagnosis-Surgery-Adjuvant} is the only hierarchical reference pathway in level 1 in subgroup 2.  Although the `Neo-adjuvant' network is discordant for patients in subgroup 2, patient pathways may still traverse it, so we include it in the model and explain its subsequent levels. However, the `Continual' network cannot appear in any patient pathway in subgroup 2 as this network is only used for subgroups 4 and 6, so we exclude it from further discussion for subgroup 2. 

\subsubsection{Level 2}

Each network in level 1 consists of subnetworks in level 2 that describe patient conditions that affect the pathway (see the second column of Table \ref{tab:activities}). If patient conditions do not affect the pathways of a level 1 network, the network will include only one subnetwork in level 2 with the same name and label. For example, `Diagnosis' and `Neo-adjuvant' have only one subnetwork. All subnetworks are listed by their desired chronological order. Similar to level 1, these subnetworks are labeled mandatory, optional or discordant (these labels are indicated in parentheses in the table), but we also add a fourth label `alternative mandatory', which is given to subnetworks that are mutually exclusive and collectively exhaustive such that a patient must traverse one (and only one) of them. For example, the first four subnetworks of the `Surgery' network are `alternative mandatory' because the surgery can be either a BCS or mastectomy, and the margins can be either positive or negative.

We generate hierarchical reference pathways in level 2 using a combination of level 2 subnetworks and their labels. For example, 
the set of hierarchical reference pathways for the level 1 `Surgery' network are \textit{BCS with neg margin - Post-op}; \textit{BCS with pos margin - Post-op}; \textit{Mastectomy with neg margin - Post-op}; and \textit{Mastectomy with pos margin - Post-op}.

\subsubsection{Level 3}
Networks in level 2 contain subnetworks in level 3, which describe the healthcare encounters. The third column of Table \ref{tab:activities} shows the subnetworks of level 2 networks and their requirements. 
Networks in this level also have `mandatory', `optional', and `discordant' labels depending on their requirements, but their labels depend on the activities that they contain. We describe the details of the activities in Section \ref{sec:activities}. 

Similar to the previous levels, we generate hierarchical reference pathways in level 3 using the labels of level 3 subnetworks. For example, the hierarchical reference pathways inside the level 2 `Diagnosis' network are \textit{Consultation- Screening Imaging- Diagnostic Imaging- Tissue Diagnosis- Pre-op Imaging}; \textit{Consultation- Diagnostic Imaging- Tissue Diagnosis- Pre-op Imaging}; \textit{Consultation- Screening Imaging- Diagnostic Imaging- Tissue Diagnosis}; and \textit{Consultation- Diagnostic Imaging- Tissue Diagnosis}.

\subsubsection{Level 4}
\label{sec:activities} 
Networks in level 3 do not include any subnetworks, but consist of clinical activities that are in level 4. 
The last column of Table \ref{tab:activities} presents the complete set of activities for each network in level 3 and the pathway requirements. 

Patients with high breast cancer risk are advised to undergo screening mammography every year, but it is not mandatory, so we label it optional. 
All patients in the cohort are symptomatic or an abnormality was observed in their screening mammogram, so they receive a mandatory diagnostic mammogram. A breast ultrasound and/or ductogram are optional. Sampling of the patient's breast tissue is also needed to confirm the diagnosis, so a breast biopsy is mandatory. In some cases, lymphatic node biopsies might be needed during the diagnostic evaluation and before surgery, so they are labeled as optional. A surgeon consultation is also required to prepare the patient for the upcoming surgery, so it is mandatory. The consultation can happen any time before the surgery. Shortly before undergoing surgery, patients may have an optional chest imaging procedure to identify comorbidities that may affect the conduct of anesthesia for the surgery. We use these requirements to label the activities of level 3 subnetworks that are subnetworks of `Diagnosis' network in level 2. 

On the day of the surgery, we have a mandatory activity for surgery preparation but optional activities for breast ultrasound and mammogram for patients who have `Imaging to guide surgery'. The surgery is either a mastectomy or BCS; if we see either we assume it is a concordant activity. A mandatory SLNB must be performed in which the sentinel lymph node is identified, removed, and examined. This operation helps with cancer staging and guiding the remainder of the treatment. If a patient has positive margins after the surgery, meaning cancer cells are seen at the edge of the removed tissue, a re-excision is required. This (mandatory) follow-up surgery is necessary and is not considered discordant. After surgery, patients may receive an optional bone scan to determine whether cancer has spread to the bones. Before adjuvant therapy begins, a consultation with a medical oncologist and a radiation oncologist is mandatory. We use these requirements to label the activities in level 3 subnetworks that are subnetworks of `BCS with neg margin', `BCS with pos margin', `Mastectomy with neg margin', `Mastectomy with pos margin', and `Post-op' in level 2. 

Since the patients in this subgroup are ER-positive, hormone therapy is mandatory, but chemotherapy is optional. Due to data limitations, we cannot track all hormone therapy events, so we overestimate concordance and assume that the patient is concordant if we observe a single hormone therapy event in the patient's pathway. For chemotherapy, at least four cycles is considered concordant (i.e., chemotherapy was deemed to have been completed).  Radiation therapy is used to kill any cancer cells that may remain after surgery if the patient had BCS or if they had positive margins, and thus is mandatory for these patients. Determining concordance for radiation therapy depends on the number of treatment sessions, which depends on the dose prescription. Our data did not record the dose prescription, so we assumed that a patient who had 12-16 treatment events (high dose per session) or greater than 20 treatment events (low dose per session) was concordant (i.e., radiation was deemed to have been completed). Any other number was considered discordant. We use these requirements to label the activities that are in level 3 subnetworks that are subnetworks of `BCS or pos margin' and `Mastectomy and neg margin' networks in level 2. Note that we create two activity nodes for each therapy activity to indicate whether therapy is completed or not. In addition, chemotherapy is optional for patients, so its corresponding network in level 3 is optional, but if the patient pathway traverses this network, it needs to have chemotherapy activities, making them mandatory.

We generate hierarchical reference pathways in level 4 using the combinations of the activities and their labels as mentioned before. For example, \textit{Surgery consultation} is the only hierarchical reference pathway for the level 3 `Consultation' subnetwork of the level 2 `Diagnosis' network. 

The procedure of modeling the requirements of all the other subgroups and generating their hierarchical reference pathways is similar. The requirement tables for all subgroups are available in Appendix \ref{app:tables}. 

\onecolumn
\begin{ThreePartTable}
\begin{longtable}{llll}
\caption{Clinical pathway requirements for subgroup 2} \label{tab:activities}\\

\toprule
Level 1          & Level 2                        & Level 3                & Level 4                   \\
\toprule

Diagnosis (M)    & Diagnosis (M)       & Consultation (M)                  & Surgery Consultation (M)         \\
                 &                     &                              & Chemotherapy Consultation (D)              \\
                 &                     &                              & Radiation Consultation (D)              \\
                \cmidrule(l){3-4}

                 &                     & Screening Imaging (O)        & Screening mammogram (O)     \\
                 \cmidrule(l){3-4}
                 &                     & Diagnostic Imaging (M)       & Breast Ultrasound (O)       \\
                 &                     &                              & Ductogram (O)               \\
                 &                     &                              & Mammogram (M)               \\
                 \cmidrule(l){3-4}
                 &                     & Tissue Diagnosis (M)                & Breast Biopsy (M)           \\
                 &                     &                              & Lymphatic Node Biopsy (O)   \\
                 \cmidrule(l){3-4}
                 &                     & Staging (D)                  & Bone Scan (D)               \\
                 &                     &                              & Abdomen CT/US (D)           \\
                 &                     &                              & Chest CT/X-ray (D)          \\
                 \cmidrule(l){3-4}
                 &                     & Pre-op Imaging (O)           & Chest X-ray (O)             \\
                 \midrule
Neo-adjuvant (D) & Neo-adjuvant (D)    & Hormone therapy (D)          & Hormone therapy-start (D)         \\
                 &                     &                              & Hormone therapy-comp (D)    \\
                 \cmidrule(l){3-4}
                 &                     & Chemotherapy (D)             & Chemotherapy-start (D)            \\
                 &                     &                              & Chemotherapy-comp (D)       \\
                 \cmidrule(l){3-4}
                 &                     & Targeted therapy (D)         & Targeted therapy-start (D)        \\
                 &                     &                              & Targeted therapy-comp (D)   \\
                 \midrule
Surgery (M)      & BCS with            & Imaging to guide surgery (O) & Preparation for surgery (M) \\
                 & neg margin (AM)     &                              & Breast Ultrasound (O)       \\
                 &                     &                              & Mammogram (O)               \\
                 \cmidrule(l){3-4}
                 &                     & Surgery (M)                  & BCS (M)                     \\
                 \cmidrule(l){3-4}
                 &                     & Management of Axilla (M)          & SLNB (M)                    \\
                 &                     &                              & ALND (D)                    \\
                 \cmidrule(l){2-4}
                 & BCS with            & Imaging to guide surgery (O) & Preparation for surgery (M) \\
                 & pos margin (AM)     &                              & Breast Ultrasound (O)       \\
                 &                     &                              & Mammogram (O)               \\
                 \cmidrule(l){3-4}
                 &                     & Surgery (M)                  & BCS (M)                     \\
                 \cmidrule(l){3-4}
                 &                     & Management of Axilla (M)          & SLNB (M)                    \\
                 &                     &                              & ALND (D)                    \\
                 \cmidrule(l){3-4}
                 &                     & Follow-up surgery (M)        & Follow-up surgery (M)       \\
                 \cmidrule(l){2-4}
                 & Mastectomy with     & Imaging to guide surgery (O) & Preparation for surgery (M) \\
                 & neg margin (AM)     &                              & Breast Ultrasound (O)       \\
                 &                     &                              & Mammogram (O)               \\
                 \cmidrule(l){3-4}
                 &                     & Surgery (M)                  & Mastectomy (M)              \\
                 \cmidrule(l){3-4}
                 &                     & Management of Axilla (M)          & SLNB (M)                    \\
                 &                     &                              & ALND (D)                    \\
                                  \cmidrule(l){2-4}

                 & Mastectomy with     & Imaging to guide surgery (O) & Preparation for surgery (M) \\
                 & pos margin (AM)     &                              & Breast Ultrasound (O)       \\
                 &                     &                              & Mammogram (O)               \\
                 \cmidrule(l){3-4}
                 &                     & Surgery (M)                  & Mastectomy (M)              \\
                 \cmidrule(l){3-4}
                 &                     & Management of Axilla (M)          & SLNB (M)                    \\
                 &                     &                              & ALND (D)                    \\
                 \cmidrule(l){3-4}
                 &                     & Follow-up surgery (M)        & Follow-up surgery (M)       \\
                    \cmidrule(l){2-4}

                 & Post-op (M)         & Consultation (M)                  & Medical Oncologist Consultation (M)              \\
                 &                     &                              & Radiation Oncologist Consultation (M)              \\
                 \cmidrule(l){3-4}
                 &                     & Staging (O)                  & Bone Scan (O)               \\
                 \midrule
Adjuvant (M)     & BCS or              & Hormone therapy (M)          & Hormone therapy-start (M)         \\
                 & pos margin (AM)     &                              & Hormone therapy-comp (M)    \\
                 \cmidrule(l){3-4}
                 &                     & Chemotherapy (O)             & Chemotherapy-start (M)            \\
                 &                     &                              & Chemotherapy-comp (M)       \\
                 \cmidrule(l){3-4}
                 &                     & Targeted therapy (D)         & Targeted therapy-start (D)        \\
                 &                     &                              & Targeted therapy-comp (D)   \\
                 \cmidrule(l){3-4}
                 &                     & Radiation (M)                & Radiation-start (M)               \\
                 &                     &                              & Radiation-comp (M)           \\
                    \cmidrule(l){2-4}

                 & Mastectomy and         & Hormone therapy (M)          & Hormone therapy-start (M)         \\
                 & neg margin (AM) &                              & Hormone therapy-comp (M)    \\
                 \cmidrule(l){3-4}
                 &                     & Chemotherapy (O)             & Chemotherapy-start (M)            \\
                 &                     &                              & Chemotherapy-comp (M)       \\
                 \cmidrule(l){3-4}
                 &                     & Targeted therapy (D)         & Targeted therapy-start (D)        \\
                 &                     &                              & Targeted therapy-comp (D)   \\
                 \cmidrule(l){3-4}
                 &                     & Radiation (D)                & Radiation-start (D)               \\
                 &                     &                              & Radiation-comp (D)        \\
\bottomrule
\end{longtable}
\begin{tablenotes}
        \scriptsize
        \item \textbf{M}: mandatory, \textbf{O}: optional, \textbf{D}: discordant, \textbf{AM}: alternative mandatory
    \end{tablenotes}
\end{ThreePartTable}
\twocolumn

\subsection{Summary of methodology}
\label{sec:sum}
Before proceeding to our numerical results, we provide a brief summary of the overall methodology. We first form the hierarchical network using Table \ref{tab:activities} or Tables \ref{tab:activities3}-\ref{tab:activities9}. The hierarchical network structure is the same for all patient subgroups. Then, for the included subgroups, we enumerate reference pathways for subnetworks across all levels. A reference pathway includes all mandatory subnetworks and permutations of optional ones in the order that appear in the tables. We refine patient pathways and include activities between 6 months before diagnosis until 15 months after diagnosis. We then map patient and hierarchical reference pathways to the network and represent them as flow vectors. We divide patient pathways into training (diagnosed in 2010-2012) and testing sets (diagnosed in 2013-2016). After solving model~\eqref{model:hierioref} with the reference pathways, the training set is used as input to model~\eqref{model:hierpatientref} to find the final arc costs, which form the concordance metric. The testing set is used in the survival and points of discordance analyses to evaluate the concordance metric out of sample.

The following section explains how we use the reference and patient pathways in the inverse optimization problems.

\section{Results}
\label{sec:results}

\subsection{Inverse Optimization Results}
\label{sec:ioresults}

\begin{figure*}
\centering  
\subfigure[Subgroup 2]{\includegraphics[width=0.4\linewidth]{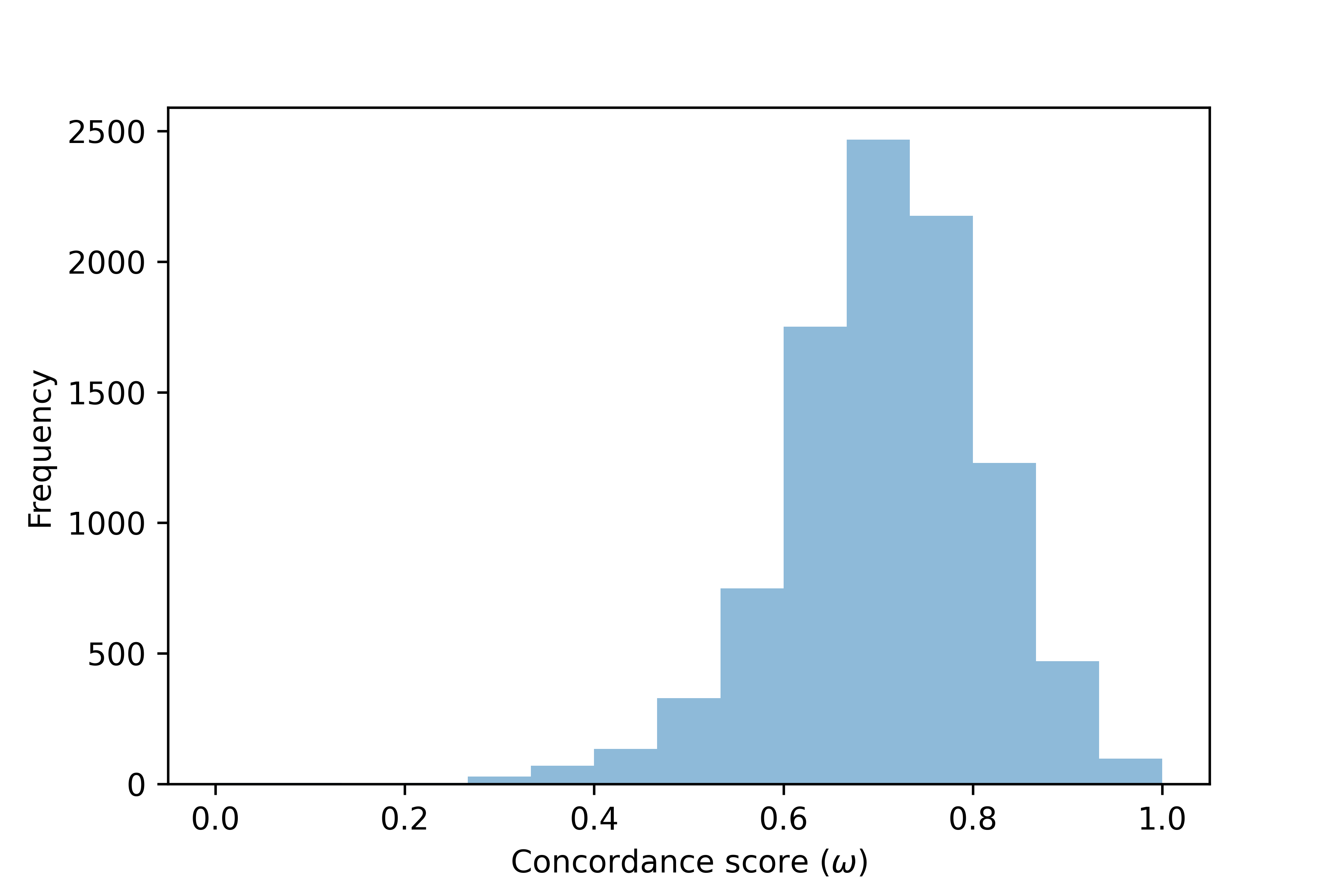}}
\subfigure[Subgroup 3]{\includegraphics[width=0.4\linewidth]{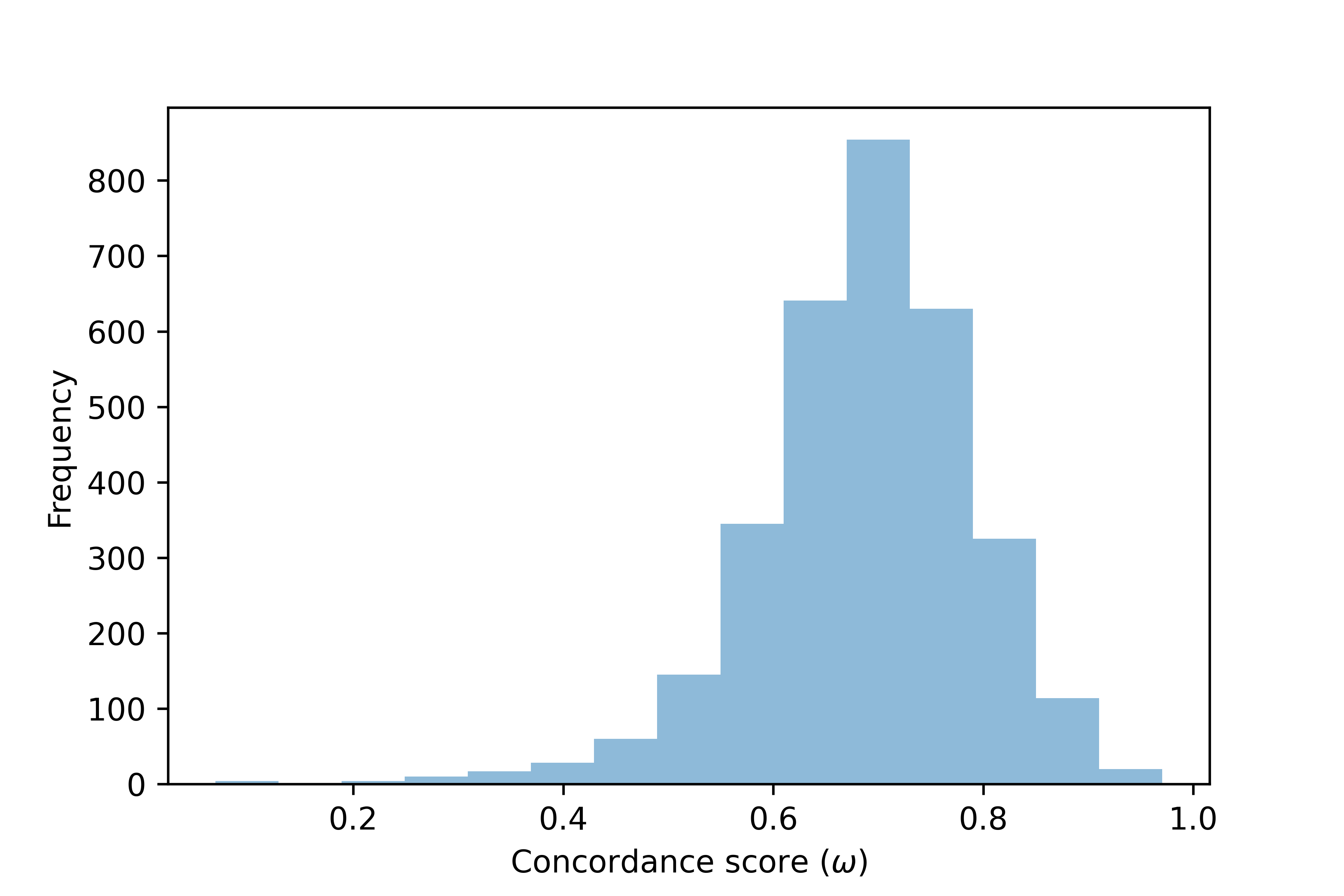}}
\subfigure[Subgroup 4]{\includegraphics[width=0.4\linewidth]{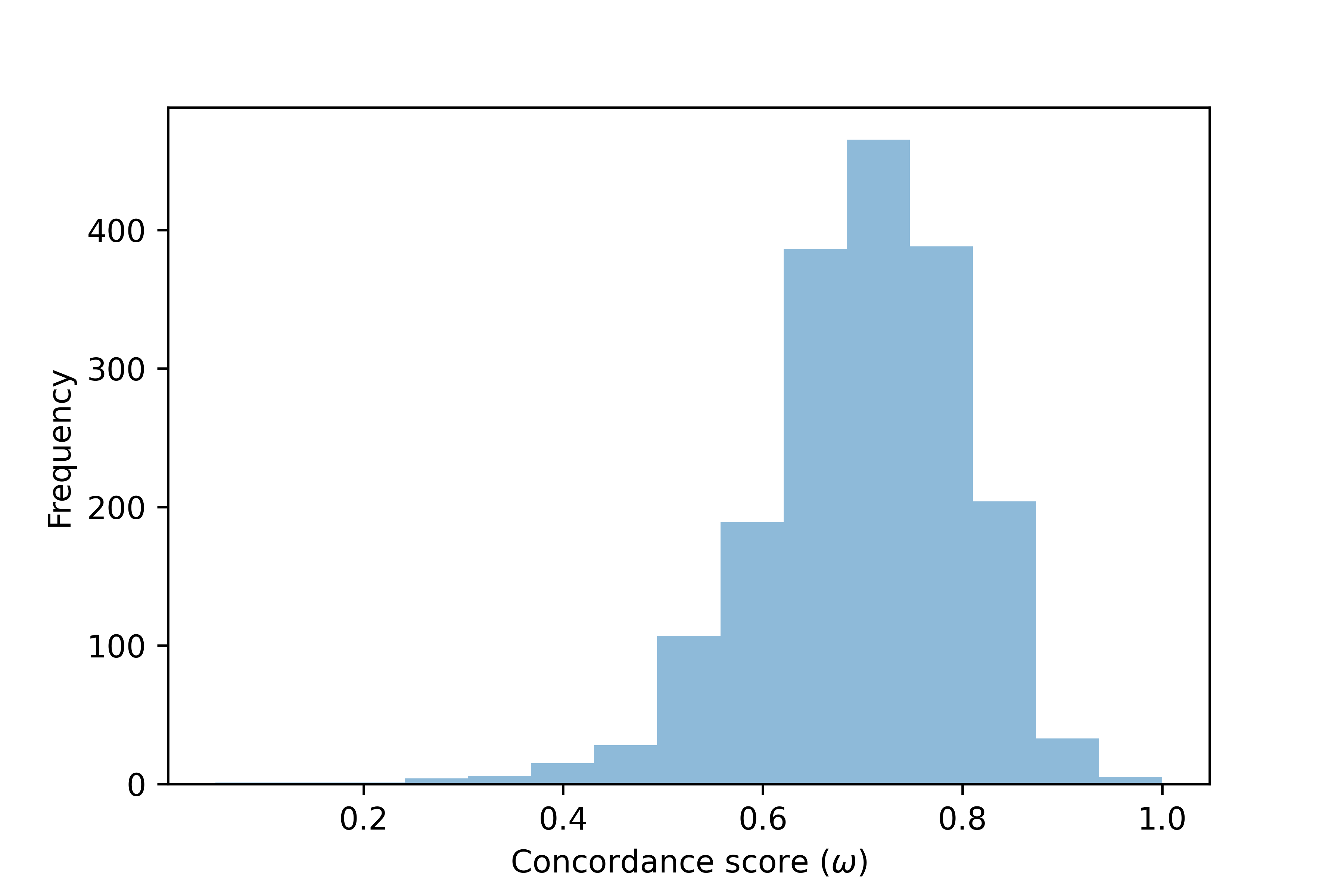}}
\subfigure[Subgroup 5]{\includegraphics[width=0.4\linewidth]{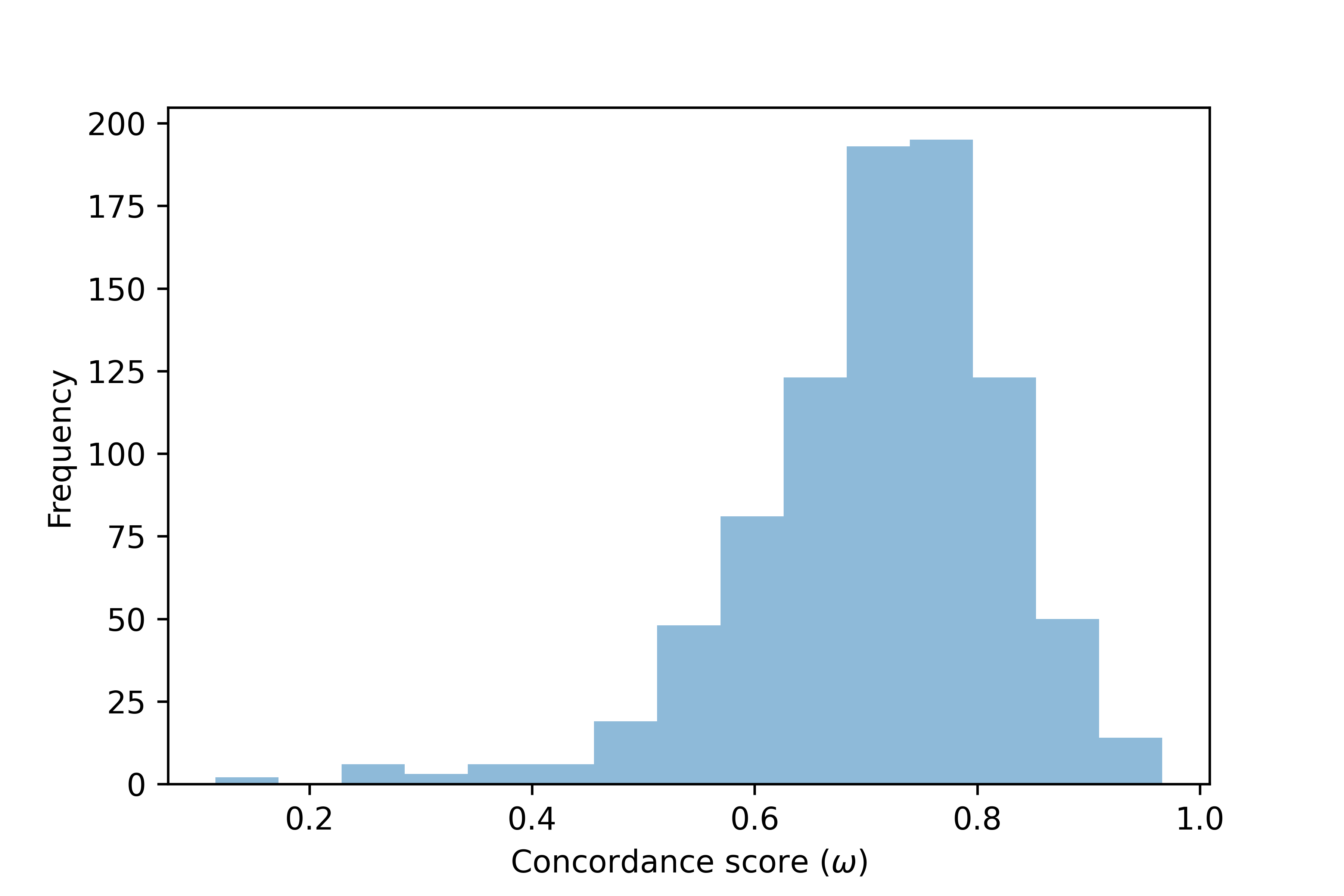}}
\subfigure[Subgroup 6]{\includegraphics[width=0.4\linewidth]{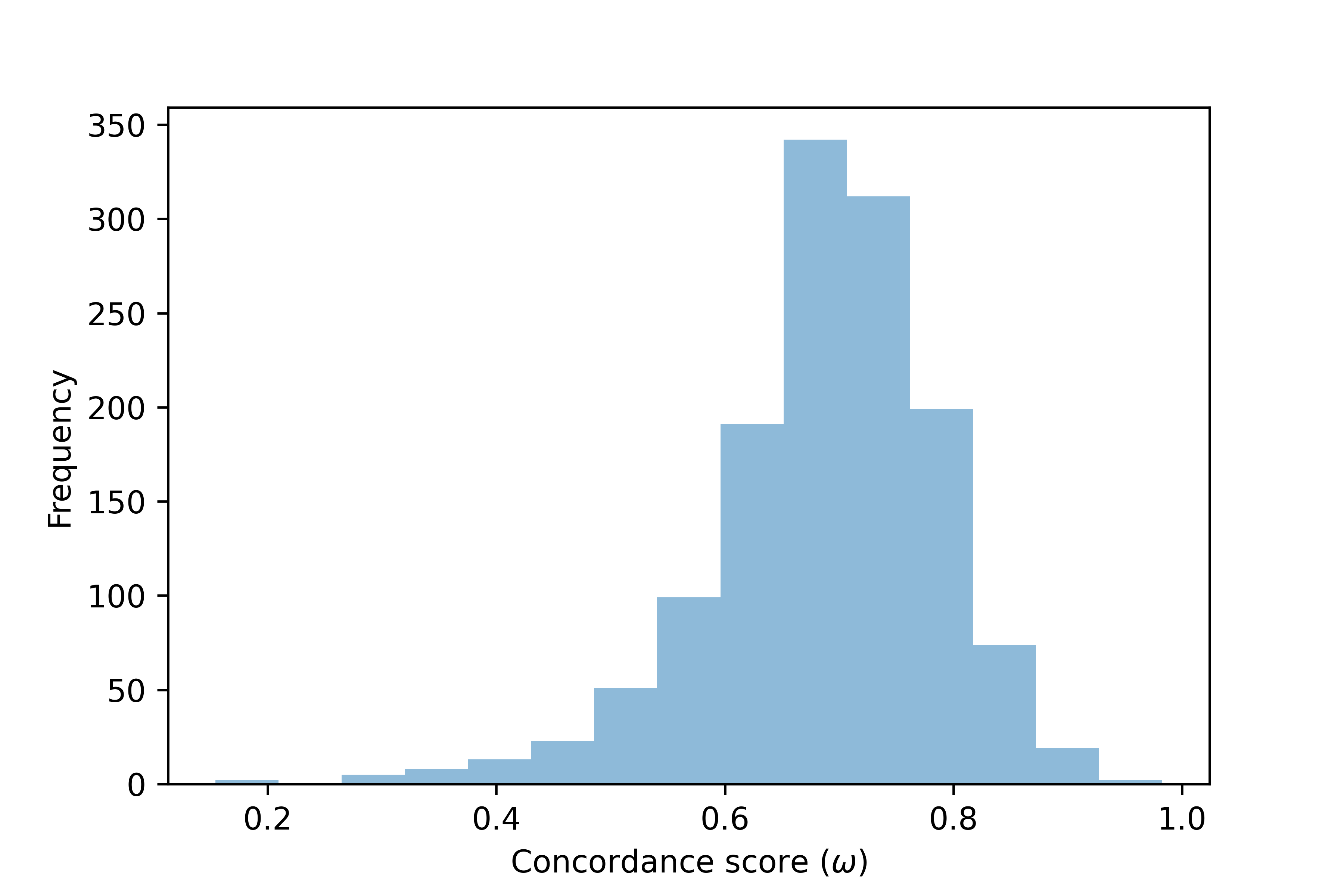}}
\subfigure[Subgroup 7]{\includegraphics[width=0.4\linewidth]{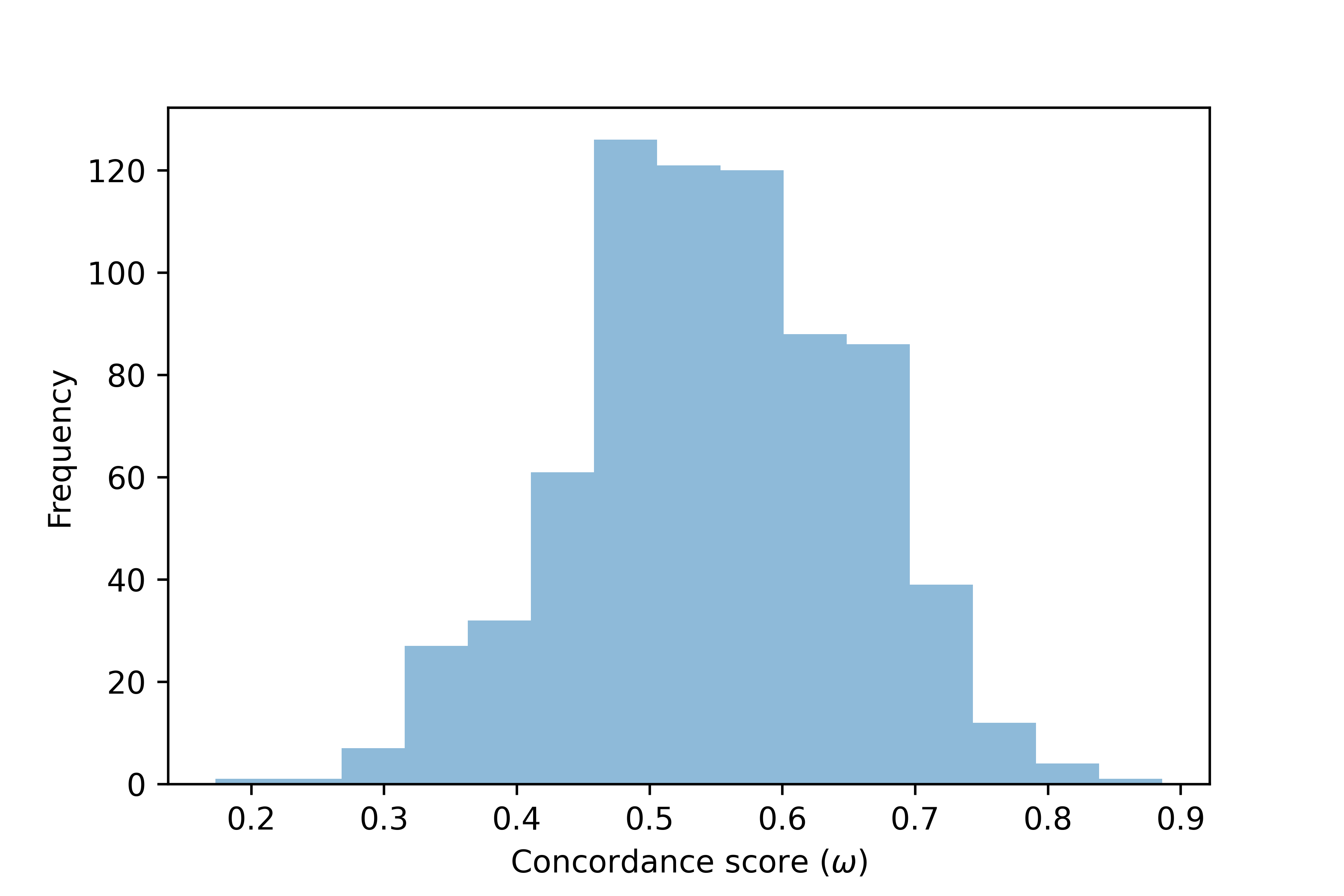}}
\subfigure[Subgroup 8]{\includegraphics[width=0.4\linewidth]{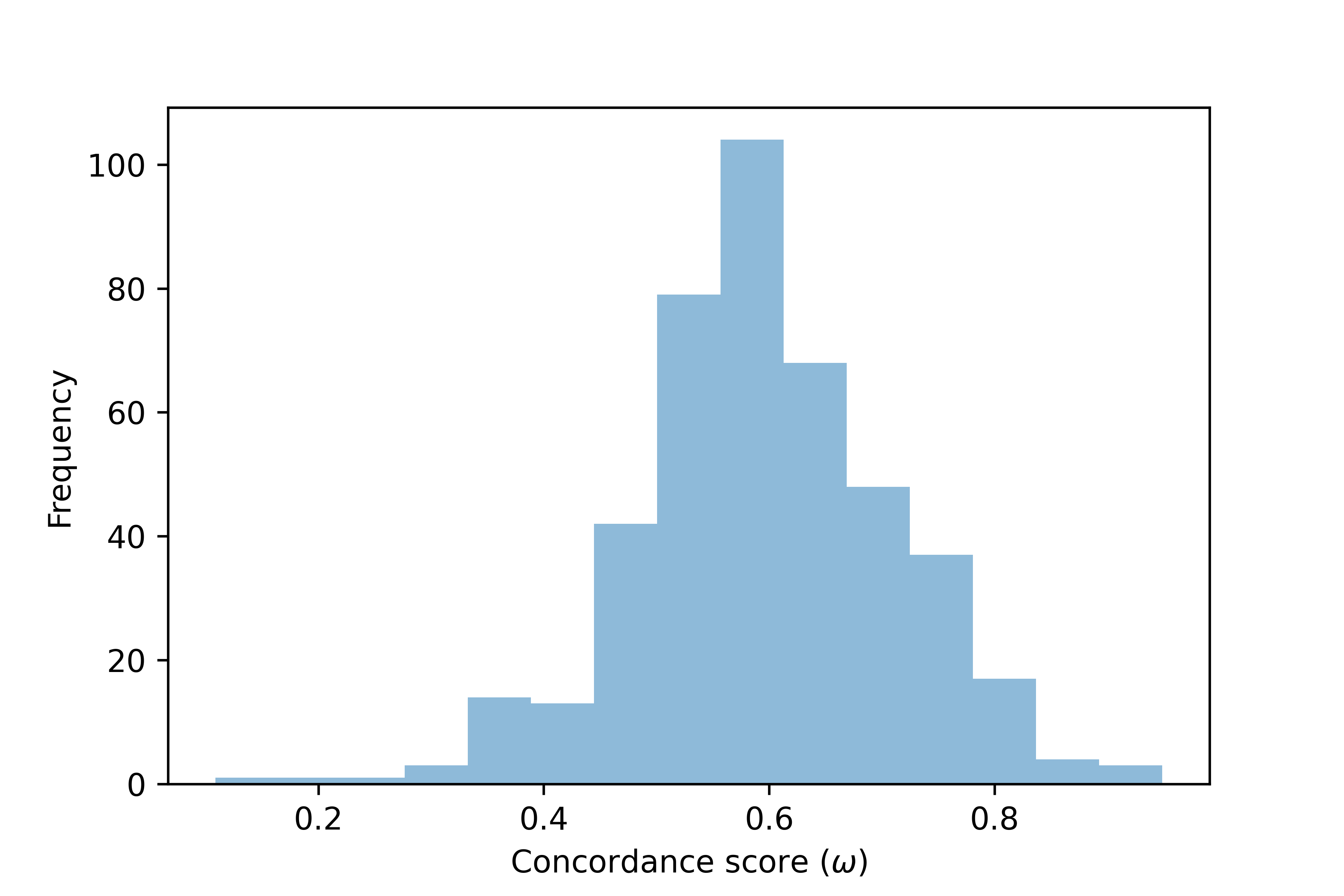}}
\subfigure[Subgroup 9]{\includegraphics[width=0.4\linewidth]{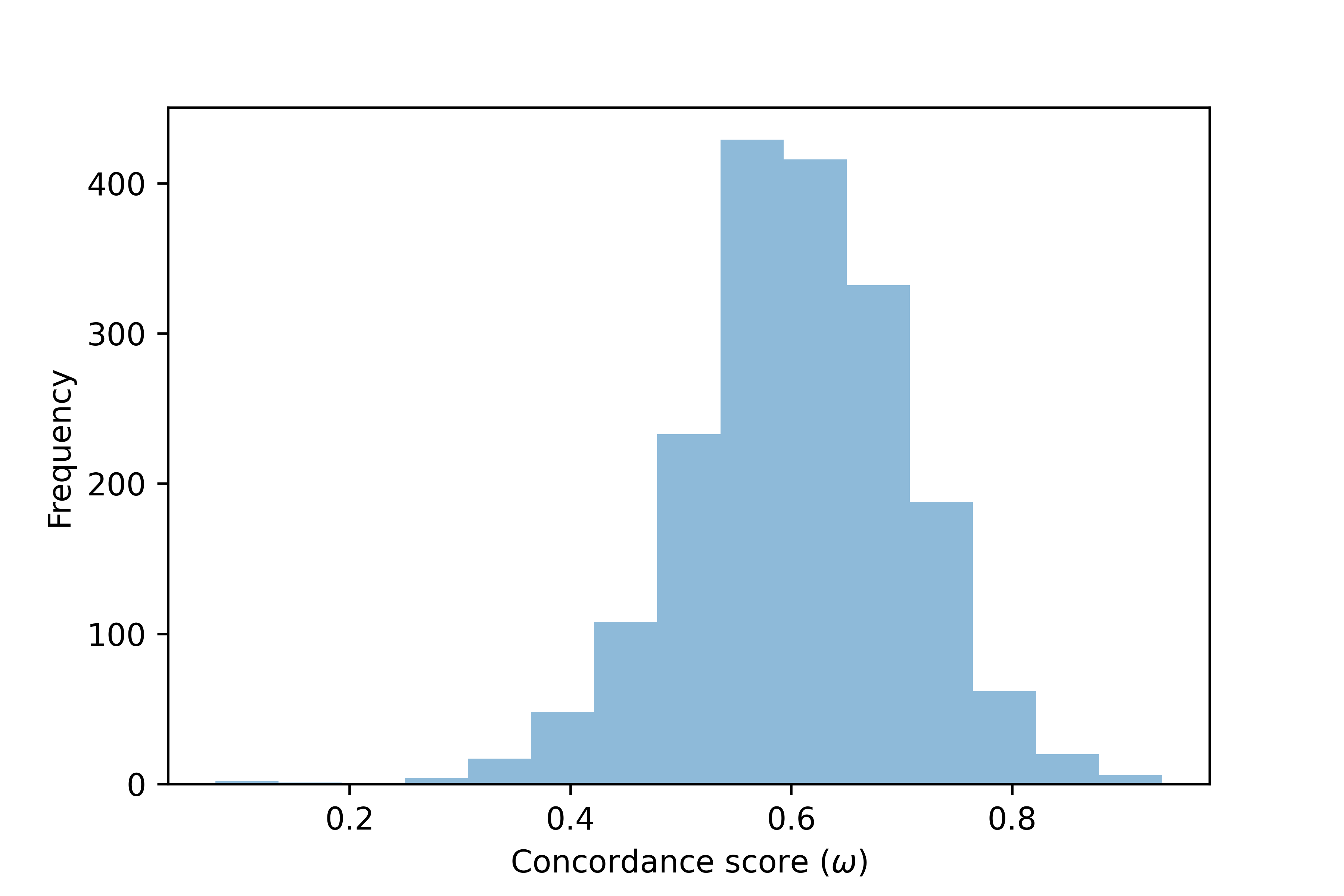}}
\caption{Distribution of concordance score in training sets}
\label{fig:dists}
\end{figure*}

The hierarchical network, hierarchical reference pathways, and patient pathways from the training dataset are the inputs to inverse optimization formulations \eqref{model:hierioref} and \eqref{model:hierpatientref}. For each subgroup, we solve formulation \eqref{model:hierioref} with the hierarchical reference pathways, followed by formulation \eqref{model:hierpatientref} with the survived and died patient data, to obtain an optimal cost vector. 
Then, we calculate concordance scores ($\omega$) for all patient pathways in each subgroup using equation \eqref{eq:omega}. Figure \ref{fig:dists} shows the distribution of the concordance scores in the training dataset. 
A concordance score of 1 indicates that the patient pathway is a shortest path and has followed the clinical pathway. According to the figure, almost no patient has followed any of the clinical pathways.
We observe that concordance score distributions of subgroups 2-6 are left skewed (mode between 0.6 and 0.8) while those of subgroups 7-9 are approximately normal (mode between 0.5 and 0.6). \ny{This indicates that a larger proportion of patients have high concordance scores in subgroups 2-6.} 

To test the robustness of our results to the training data, we repeatedly solved model \eqref{model:hierpatientref} with bootstrapped samples of the training data, and the resulting concordance scores were essentially unchanged \ny{(see Figure \ref{fig:boot} in Appendix \ref{app:fig_table})}.

All computations are conducted in Python with Gurobi version 9.0.3 to solve the optimization problems on a Macbook with Dual-Core Intel Core i7 and 16 GB of RAM.  The code is available in our Github repository at \url{https://github.com/yusufshalaby/DPM-Breast-Package}.

\subsection{Survival Analysis}
\label{sec:survival}

In this section, we demonstrate that there is a statistically significant association between survival and concordance for each subgroup using the testing data (2013-2016). This result provides an out-of-sample validation of the concordance metric.

We used a Cox proportional hazards model \citep{cox} to evaluate the association between survival and concordance, considering potential confounders and predictors of survival based on standard practice or the evidence that the variable is correlated with both concordance and survival. These variables consist of sociodemographic factors (age, sex, rural residency, quintiles of neighbourhood income, terciles of neighbourhood immigrant population density), Charlson comorbidity score \citep{charlson1987new,quan2005coding}, cancer diagnostic and treatment characteristics (screening group), mental health history (psychotic and non-psychotic disorders, substance use disorders, other social, family, and occupational issues) from three years before diagnosis until one year after diagnosis, and healthcare utilization (number of OHIP billing dates for one year before six months of diagnosis, number of hospital admissions six months before diagnosis, number of unscheduled ED visits six months before diagnosis), which is a proxy for overall patient health. 

Based on a preliminary variable selection \cite{ishwaran2008random} and checking the Cox model assumptions, we used age (in years), rural residency (rural, urban), Charlson score (0, 1, 2+), neighbourhood income (highest, mid-high, mid-low, middle, lowest), neighbourhood immigrant population density (most dense, mid dense, least dense), screening (screened by Ontario Breast Screening Program \cite{obsp} (OBSP-screened), screened by a general practitioner (GP-screened), symptomatic), and two mental health indicators (psychotic disorders, substance use disorders) in the model and examined four-year survival. We excluded patients who did not have all these covariate values. We also excluded patients who died within 15 months of diagnosis because they may not have had the opportunity to complete their treatment. Table \ref{tab:datasurvival} summarizes the characteristics of the cohort used in the survival analysis for each subgroup.

\begin{table*}
\scriptsize
\caption{Characteristics of cohort included in survival analysis}
\label{tab:datasurvival}
\begin{tabular}{L{3.4cm}llllllll}
\toprule
& \multicolumn{8}{c}{Subgroup}       \\
                                 \cmidrule(r){2-9}                                                                   &\ 2  \ & \ 3  \ & \ 4  \ & \ 5  \ & \ 6  \ & \ 7  \ & \ 8  \  & \ 9  \  \\
                        \toprule
Total                             & 13,718 & 4,403 & 2,520 & 1,313 & 2,198 & 971  & 501  & 2,417 \\
\midrule
4-year mortality, n (\%)                             & 552 (4.0)   & 219 (5.0)  & 253 (10.0)  & 29 (2.2)   & 145 (6.6) & 144 (14.8)  & 211 (42.1) & 351 (14.5)  \\
Concordance score, mean (SD)                 & 0.74 (0.1) & 0.70 (0.1) & 0.72 (0.1)& 0.74 (0.1)& 0.71 (0.1)& 0.58 (0.1)& 0.62 (0.1)& 0.62 (0.1) \\
\midrule
Age, mean (SD)                               & 63.8 (12.5) & 60.3 (13.2) & 59 (13.9)   & 59.7 (11.8) & 57.5 (13.6) & 56.8 (14.1) & 57.6 (14.2) & 60.2 (14.6) \\
Female, n (\%)                            & 13,632 (99.4) & 4,354 (98.9) & 2,517 (99.9) & 1,310 (99.8) & 2,187 (99.5) & 965 (99.4) & 501 (100)  & 2,393 (99.0) \\
Rural residency, n (\%)                   & 1,621 (11.8) & 489 (11.1) & 289 (11.5) & 141 (10.7) & 232 (10.6) & 107 (11.0) & 52 (10.4)  & 267 (11.0) \\
Charlson score=0, n (\%)                  & 12,114 (88.3) & 3,904 (88.7) & 2,216 (87.9) & 1,197 (91.2) & 1,979 (90.0) & 871 (89.7) & 449 (89.6)  & 2,120 (87.7) \\
\midrule
\multicolumn{9}{l}{Neighbourhood Income, n (\%)}       \\
\quad Highest                           & 3,235 (23.6)  & 973 (22.1)  & 542 (21.5)  & 301 (22.9)  & 481 (21.9)  & 206 (21.2)  & 103 (20.6)  & 537 (22.2) \\
\quad Mid-high                          & 2,991 (21.8)  & 971 (22.1)  & 577 (22.9)  & 318 (24.2)  & 498 (22.7)  & 211 (21.7)  & 107 (21.4)  & 524 (21.7)  \\
\quad Mid-low                           & 2,538 (18.5)  & 850 (19.3)  & 485 (19.2)  & 246 (18.7)  & 417 (19.0)  & 185 (19.1)  & 93 (18.6)   & 441 (18.2)  \\
\quad Middle                            & 2,746 (20.0)  & 865 (19.6)  & 511 (20.3)  & 260 (19.8)  & 446 (20.3)  & 180 (18.5)   & 107 (21.4) & 475 (19.7) \\
\quad Lowest                            & 2,208 (16.1)  & 744 (16.9)  & 405 (16.1)  & 188 (14.3)  & 356 (16.2)  & 189 (19.5)  & 91 (18.2)  & 440 (18.2)  \\
\midrule
\multicolumn{9}{l}{Neighbourhood  immigration density, n (\%)}       \\
\quad Most dense                        & 2,256 (16.4)  & 746 (16.9) & 421 (16.7) & 228 (17.4) & 409 (18.6) & 171 (17.6) & 108 (21.6) & 439 (18.2) \\
\quad Mid dense                         & 3,406 (24.8)  & 1,147 (26.1) & 598 (23.7)  & 348 (26.5)  & 580 (26.4)  & 244 (25.1)  & 111 (22.2)  & 630 (26.1) \\
\quad Least dense                       & 8,056 (58.7)  & 2,510 (57.0) & 1,501 (59.6) & 737 (56.1)  & 1,209 (55.0) & 556 (57.3)  & 282 (56.3)  & 1,348 (55.8) \\
\midrule
Screening, n (\%)                         &       &      &      &      &      &      &      &      \\
\quad OBSP-screened                     & 5,299 (38.6)  & 1,186 (26.9) & 631 (25.0)  & 508 (38.7)  & 417 (19.0)  & 125 (12.9)  & 48 (9.6)  & 328 (13.6)  \\
\quad GP-screened                       & 2,679 (19.5)  & 779 (17.7)  & 478 (19.0)  & 261 (19.9)  & 401 (18.2)  & 171 (17.6)  & 92 (18.4)  & 452 (18.7)  \\
\quad Symptomatic                       & 5,740 (41.8)  & 2,438 (55.4) & 1,411 (56.0) & 544 (41.4)  & 1,380 (62.8) & 675 (69.5)  & 361 (72.1)  & 1,637 (67.7) \\
\midrule
Psychotic disorders, n (\%)               & 351 (2.6)  & 128 (2.9)  & 63 (2.5)  & 26 (2.0)  & 41 (1.9)  & 23 (2.4)  & 20 (4.0)  & 64 (2.6)  \\
Substance use disorders, n (\%)           & 157 (1.1)   & 63 (1.4)  & 31 (1.2)  & 13 (1.0)  & 33 (1.5)  & 11 (1.1)  & 3 (0.6)   & 28 (1.2)\\
\bottomrule
\end{tabular}
\end{table*}

Next, we perform a survival analysis for each subgroup, adjusting for the confounders selected in the preliminary variable selection. Note that we restrict the survival analysis to females only since males constitute a small percentage of the population in all subgroups (less than 1\%). We also scale the concordance score ($\omega$) to represent a unit change of 0.1. Table \ref{tab:omegasurvival} summarizes the association between concordance and survival for both unadjusted and adjusted survival models where the unadjusted model contains the concordance metric only, while the adjusted model contains all the selected covariates mentioned above in addition to the concordance metric. We observe that both models find a significant association between survival and concordance in all subgroups. \ny{For example, in the adjusted model for subgroup 2, the hazard ratio (HR) is 0.70, which shows that a 0.1 increase in concordance score is associated with an estimated 30$\%$ decrease in risk of mortality in this subgroup. Note that we do not claim a causal relationship, and the results only suggest a correlation. }

\begin{table*}
\caption{Risk of mortality associated with concordance score ($\omega$) for subgroups }
\label{tab:omegasurvival}
\begin{tabular}{ccccccc}
\toprule
         & \multicolumn{3}{c}{Unadjusted}         & \multicolumn{3}{c}{Adjusted}           \\
        \cmidrule(r){2-4} \cmidrule(l){5-7}
Subgroup & HR   & 95\% CI      & p-value          & HR   & 95\% CI      & p-value          \\
\toprule
2        & 0.55 & (0.52, 0.58) & \textless{}0.0001 & 0.70 & (0.66, 0.74) & \textless{}0.0001 \\
3        & 0.63 & (0.58, 0.69) & \textless{}0.0001 & 0.69 & (0.63, 0.75) & \textless{}0.0001 \\
4        & 0.51 & (0.46, 0.57) & \textless{}0.0001 & 0.59 & (0.53, 0.66) & \textless{}0.0001 \\
5        & 0.54 & (0.45, 0.65) & \textless{}0.0001 & 0.65 & (0.52, 0.82) & \textless{}0.0001 \\
6        & 0.46 & (0.40, 0.52) & \textless{}0.0001 & 0.62 & (0.54, 0.71) & \textless{}0.0001 \\
7        & 0.70 & (0.61, 0.81) & \textless{}0.0001 & 0.89 & (0.79, 0.99) & \textless{}0.0001 \\
8        & 0.68 & (0.61, 0.76) & \textless{}0.0001 & 0.79 & (0.69, 0.90) & \textless{}0.0001 \\
9        & 0.73 & (0.66, 0.82) & \textless{}0.0001 & 0.77 & (0.69, 0.86) & \textless{}0.0001\\
\bottomrule
\end{tabular}
\end{table*}

Table \ref{tab:allsurvival} shows the effects of the covariates on survival in the adjusted model for each subgroup. A (+) or ($-$) indicates a significantly positive or negative association with mortality at the 0.05 significance level while empty cells indicate no significant association. For subgroup 2, age, Charlson score (1 or 2+), psychotic disorder, and substance use disorder have statistically significant harmful effect while OBSP and GP screening, and high or middle income have statistically significant beneficial effect on survival. Age and Charlson score are commonly found harmful covariates for all the subgroups except for subgroup 8 and screening is the common beneficial covariate. Subgroup 8 has distinct results, which are less interpretable, due to the high mortality rate in this subgroup (42.1\%). 

\ny{Similar to \citep{chan2021inverse}, we repeated the survival analysis of each subgroup with different sample sizes and observed that the results are robust and the association remains significant except for small sample sizes (e.g. sample size $< 300$). }

\ny{Next, we repeated the survival analysis using bootstrapped samples of the training data used in formulation \eqref{model:hierpatientref}. The survival analysis results were identical to the original ones.}

\begin{table*}
\caption{Effects of the covariates on the risk of mortality}
\label{tab:allsurvival}
\begin{tabular}{llcccccccc}
\toprule
                                                                                                &                         & \multicolumn{8}{c}{Subgroup}       \\
                                 \cmidrule(r){3-10}                                                               &                         &\ 2  \ & \ 3  \ & \ 4  \ & \ 5  \ & \ 6  \ & \ 7  \ & \ 8  \  & \ 9  \  \\
                        \toprule
                                                                                                & Age                     & ($+$) & ($+$)  & ($+$) & ($+$)  & ($+$) & ($+$) &   & ($+$)  \\
                                                                                                & Rural residency         &  &   &  &   &  &  &  &   \\
                                                              \midrule                          
\multirow{2}{*}{\begin{tabular}[c]{@{}l@{}}Charlson score\\ (Baseline: 0)\end{tabular}}         & 1                       & ($+$) &   &  &   &  &   &   &   \\
                                                                                                & 2+                      & ($+$) & ($+$)  & ($+$) & ($+$)  & ($+$) & ($+$) &   & ($+$)  \\ \midrule
\multirow{2}{*}{\begin{tabular}[c]{@{}l@{}}Screening\\ (Baseline: Symptomatic)\end{tabular}}    & OBSP                    & ($-$) & ($-$)  & ($-$) &   &  & ($-$)  &   & ($-$)  \\
                                                                                                & GP                      & ($-$) &   &  &   &  &   & ($-$) & ($-$)  \\ \midrule
\multirow{2}{*}{\begin{tabular}[c]{@{}l@{}}Income\\ (Baseline: Lowest)\end{tabular}}    & Highest                 & ($-$) &   &  &   &  &  &  & ($-$)  \\
                                                                                        & Mid-high                &  &   &  &   &  &   & ($+$)  & ($-$)  \\
                                                                                        & Mid-low                 &  &   &  &   &  &   &   & ($-$)  \\
                                                                                        & Middle                  & ($-$) &   &  &   &  &   &   & ($-$)  \\ \midrule
\multirow{2}{*}{\begin{tabular}[c]{@{}l@{}}Immigration \\ (Baseline: Least dense)\end{tabular}} & Most dense              &  &   &  &   &  &   &   &   \\
                                                                                                & Mid dense               &  & ($-$) &  &   &  &   &   &   \\ \midrule
Mental health                                                                           & Psychotic disorder         & ($+$) &   &  &   &  &   &  & ($+$)  \\
                                                                                        & Substance use disorder  & ($+$) &   &  &   &  &   & ($+$)  &   \\ 
                                                \bottomrule
\end{tabular}
\end{table*}

\subsection{The Value of Model~\eqref{model:hierpatientref}}

\ny{This section examines the effect of removing model~\eqref{model:hierpatientref} from our inverse optimization approach. We repeat the survival analysis using another optimal cost vector of model~\eqref{model:hierioref}, where the cost vector is generated from solving model~\eqref{model:hierpatientref} with swapped patient data sets. In other words, the cost vector is an extreme case that minimizes the duality gaps for patients who died and maximizes them for patients who survived (i.e. the reverse of model~\eqref{model:hierpatientref}).}
 
\ny{Table \ref{tab:omegasurvival-reversed} in Appendix \ref{app:fig_table} summarizes the association between concordance and survival for unadjusted and adjusted models for all subgroups. We observe that the association is insignificant in all models. These results emphasize the importance of including patient data in our approach to generate an informative and valid concordance metric.}

\subsection{Points of Discordance}
\label{sec:pod}

In this subsection, we decompose the discordance of patient pathways and identify major areas of discordance in the patient population for all subgroups. 

We use equations \eqref{eq:epsilonm} and \eqref{eq:epsilona} to calculate the discordance score corresponding to each network at each level. Then, we find the population discordance for each network by summing them up over all patient pathways in the subgroup. We first summarize the contribution of discordance from networks in level 1 as they divide patient pathways into meaningful episodes of care. 
Table \ref{tab:pod-process} shows the contribution of discordance from each section towards the total discordance for each subgroup.
We observe that the majority of discordance comes from adjuvant therapy, diagnosis, and surgery in all subgroups. For subgroups 2-6 neo-adjuvant therapy has very little contribution to discordance while for subgroups 7-9 it highly contributes to discordance. One reason for the small share of discordance from the `Neo-adjuvant' subnetwork in subgroups 2, 3, and 5 is that neo-adjuvant therapy is discordant for these subgroups and only a small number of patients go through it. In contrast, neo-adjuvant therapy is recommended for subgroups 7 to 9. Note that neo-adjuvant therapy for subgroups 4 and 6 is captured in the `Continual' subnetwork. The `Other' category in Table \ref{tab:pod-process} contains all the other discordance that originates from missing a concordant section (e.g., missing `Diagnosis' in any subgroup) or having a discordant section (e.g., `Neo-adjuvant' in subgroups 2, 3, and 5). 

\begin{table*}
\caption{Percentage (\%) of total cohort discordance from level 1 subnetworks}
\label{tab:pod-process}
\begin{tabular}{lcccccccc}
\toprule
                     & \multicolumn{8}{c}{Subgroup}             \\ \cmidrule(r){2-9}
       & 2   & 3   & 4  & 5   & 6  & 7  & 8  & 9  \\
\toprule
Diagnosis  & 32  & 28  & 23 & 28  & 27 & 20 & 22 & 25 \\
Neo-adjuvant       & 3 & 5 & 8  & 3 & 8  & 17 & 19 & 13 \\
Surgery   & 23   & 23   & 22  & 31   & 19  & 14  & 17  & 17  \\
Adjuvant      & 39  & 40  & 40 & 35  & 36 & 35 & 30 & 33 \\
Other                & 3   & 4   & 7  & 3   & 10 & 14 & 12 & 12 \\
\midrule
Total                & 100   & 100  & 100 & 100  & 100 & 100 & 100 & 100 \\
\bottomrule
\end{tabular}
\end{table*}

\begin{figure*}
    \includegraphics[width=0.8\textwidth]{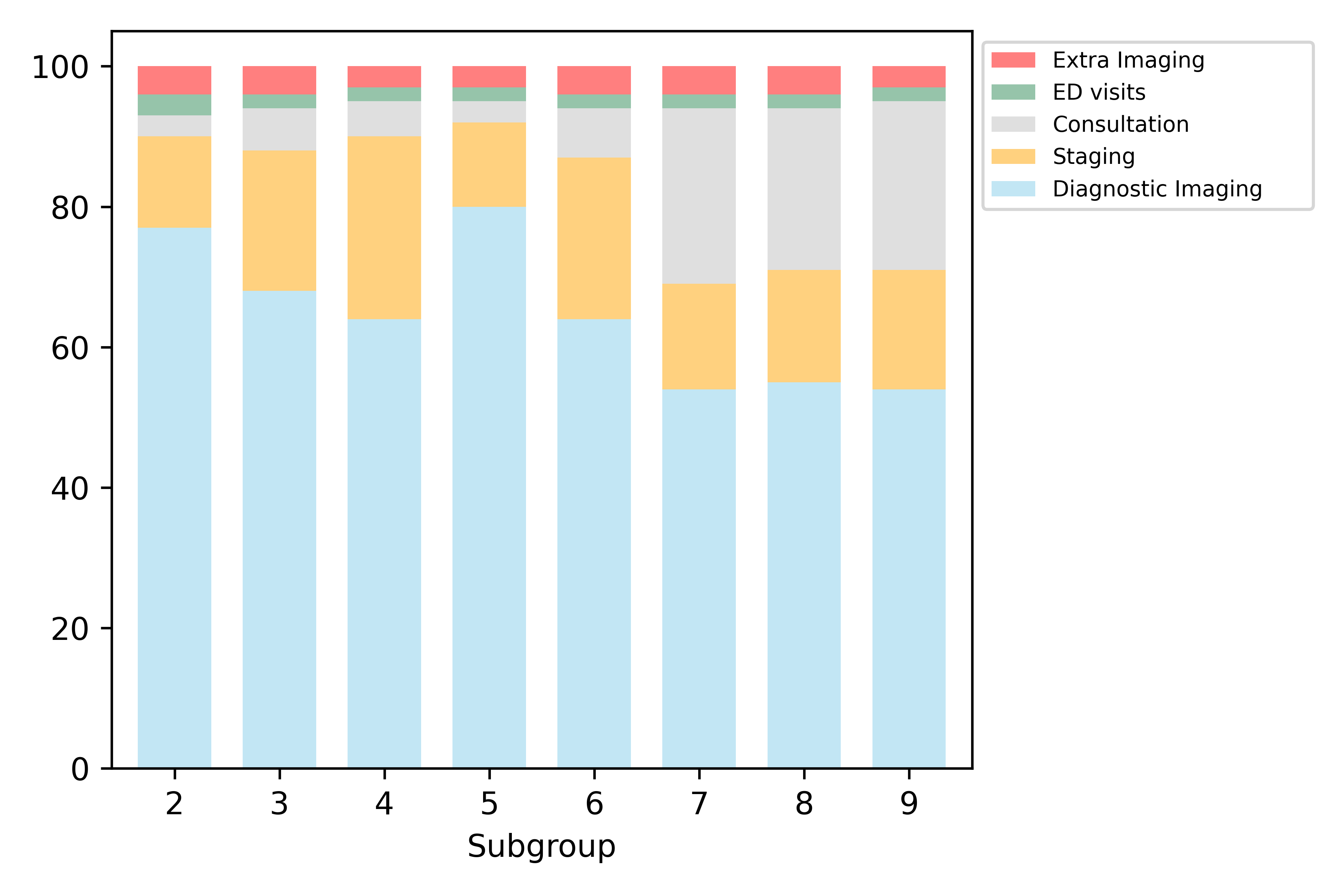}
    \caption{Frequency distribution of level 1 diagnosis subnetwork discordance over healthcare encounters (level 3 subnetworks)}
    \label{fig:diagnosis}
\end{figure*}

Next, we investigate discordance in more detail and find major points of discordance within each section. We present the distribution of discordance over healthcare encounters (i.e., networks in level 3). Figure \ref{fig:diagnosis} presents the distribution of discordance for the level 1 diagnosis subnetwork from different healthcare encounters. Healthcare encounters are either the previously introduced ones, which we listed in Table \ref{tab:activities} for subgroup 2 and in tables in Appendix \ref{app:tables} for all the other subgroups, or come from patient data and are discordant for all subgroups. A healthcare encounter is included in the figure only if its percentage of discordance is nonzero for at least one subgroup. 
Figure \ref{fig:diagnosis} shows that the main sources of discordance in the `Diagnosis' section are from `Diagnostic Imaging', `Staging', and `Consultation'. The discordance within these healthcare encounters is largely driven by extra imaging, beyond that recommended by the pathway map. `Extra Imaging' includes all the other healthcare encounters that are classified as imaging and are discordant for all subgroups.  

\begin{figure*}
    \includegraphics[width=0.8\textwidth]{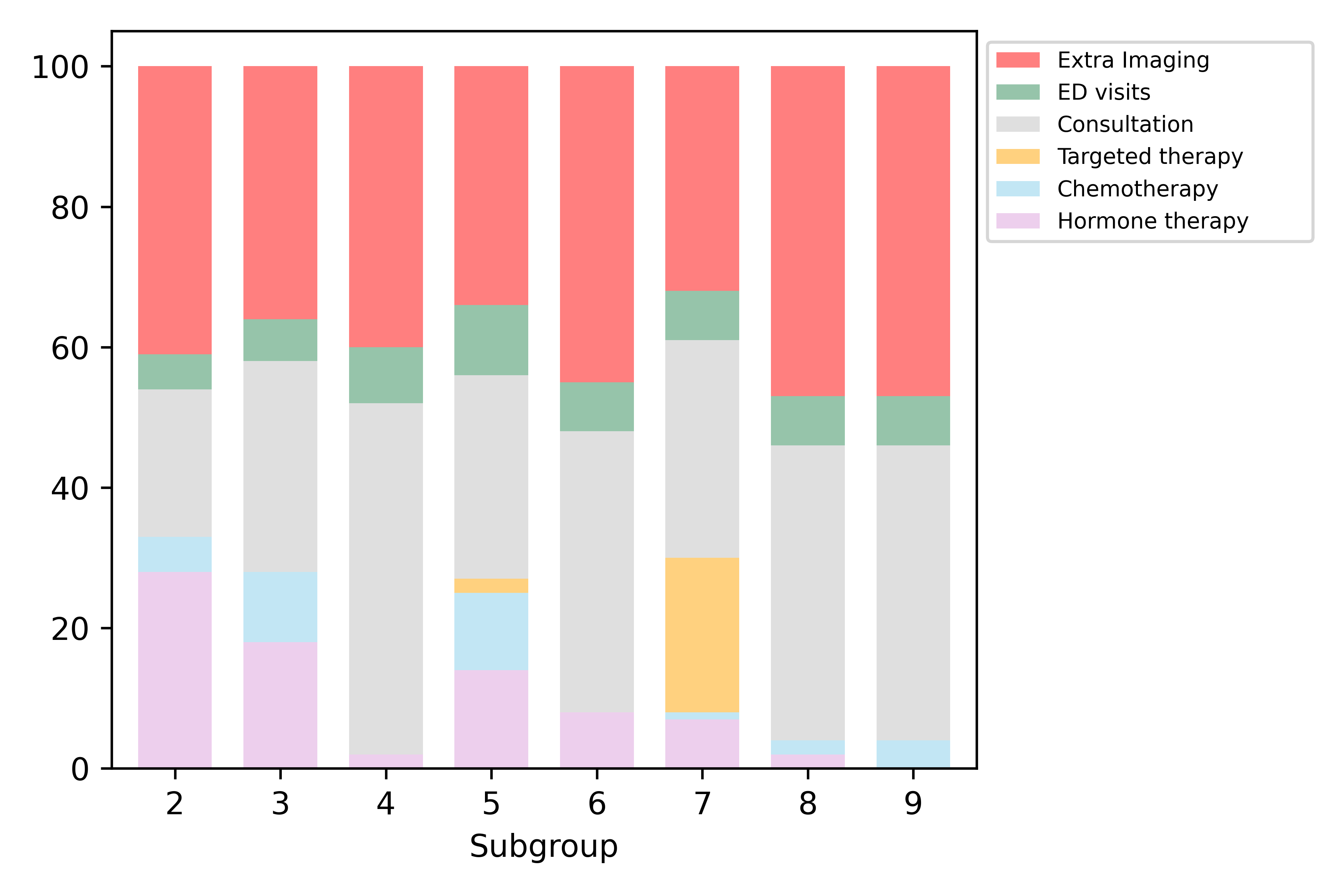}
    \caption{Frequency distribution of level 1 neo-adjuvant subnetwork discordance over healthcare encounters (level 3 subnetworks)}
    \label{fig:neo}
\end{figure*}

\begin{figure*}
    \includegraphics[width=0.8\textwidth]{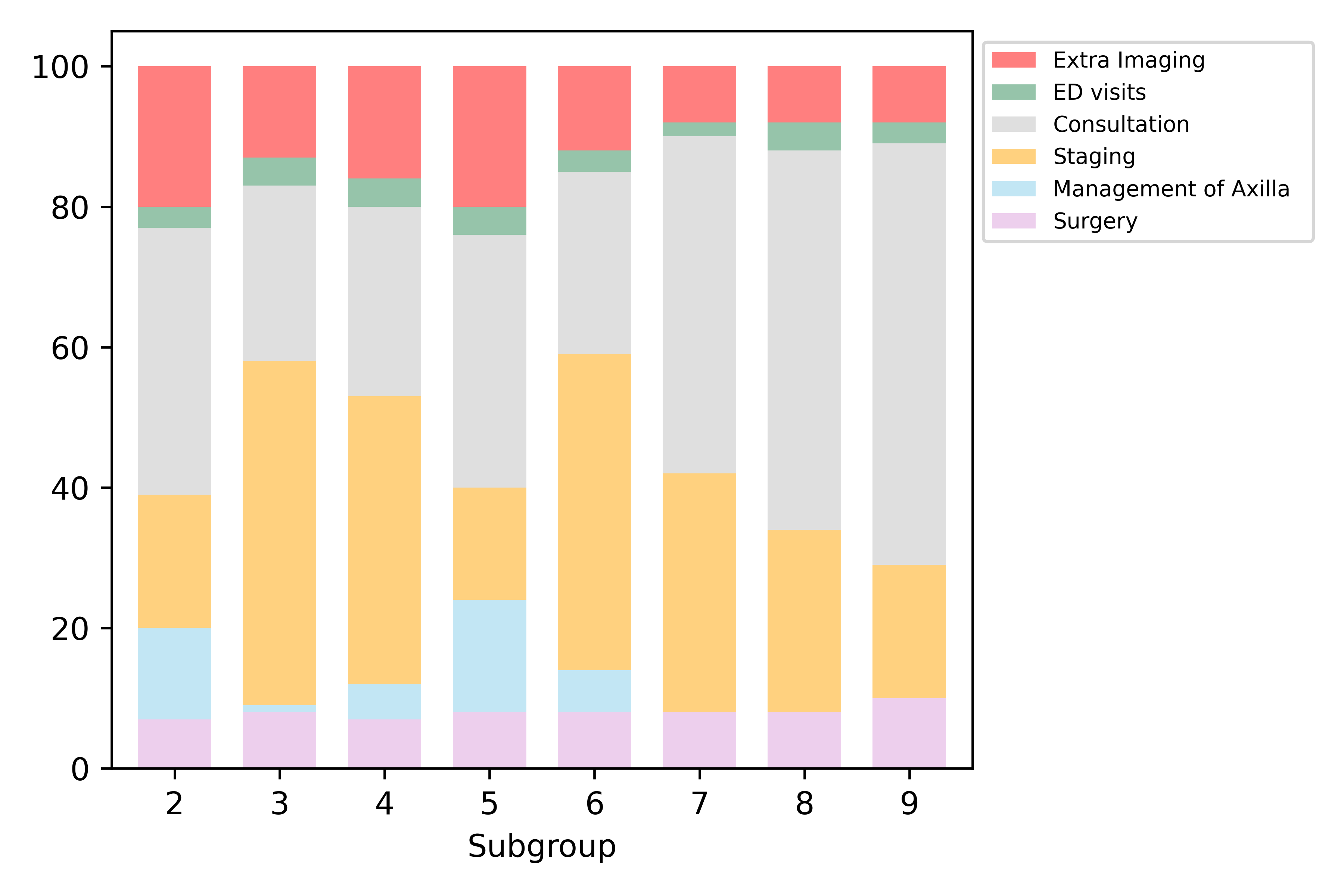}
    \caption{Frequency distribution of level 1 surgery subnetwork discordance from healthcare encounters (level 3 subnetworks)}
    \label{fig:surgery}
\end{figure*}

\begin{figure*}
    \includegraphics[width=0.8\textwidth]{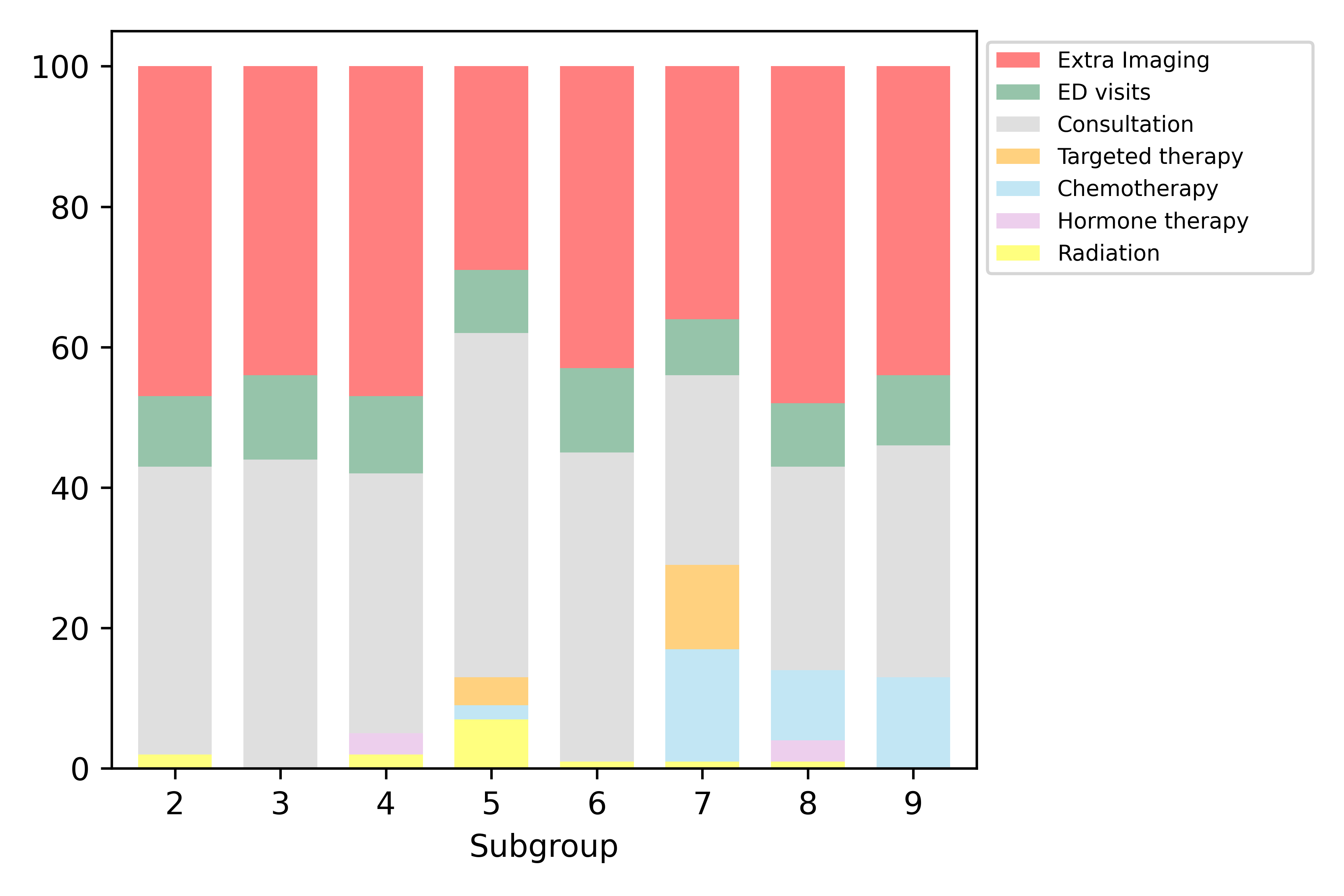}
    \caption{Frequency distribution of level 1 adjuvant subnetwork discordance from healthcare encounters (level 3 subnetworks)}
    \label{fig:adjuvant}
\end{figure*}

Figure \ref{fig:neo} presents the distribution of discordance in the level 1 `Neo-adjuvant' subnetwork. We observe that `Consultation' and `Extra Imaging' contribute the most to the discordance in this section in all subgroups, and they indicate extra consultation and imaging activities in the patient pathways. Similarly, Figure \ref{fig:surgery} suggests that `Staging', `Consultation', and `Extra Imaging' encompass the majority of discordance in the `Surgery' section in all subgroups. Finally, Figure \ref{fig:adjuvant} shows that similar to the `Neo-adjuvant' and `Surgery' subnetworks, the `Adjuvant' subnetwork includes `Consultation' and `Extra Imaging' as its main sources of discordance. Moreover, `ED visits' has a relatively high frequency in this section.

The predominance of additional activity, rather than omitted activity as primary sources of discordance is of interest. 
We conjecture that it reflects the increased comorbidity or the more numerous healthcare encounters associated with the care of those within each subgroup who die. These encounters are not articulated by the pathway map, and thus are designated as discordant. Furthermore, additional activity contributes to health system costs. Subsequent analyses of reasons for discordance and correlation with other outcomes, such as quality of life are required to identify opportunities for quality improvement. Thus, the identification of points of discordance across an entire population and the entire care continuum provides a value-based approach to prioritizing system issues for further examination and remediation.

\section{Conclusions}
\label{sec:cons}
This paper develops an inverse optimization-based concordance metric to measure clinical pathway concordance for complex diseases, where the patient journey through the healthcare system is modeled using a hierarchical network. We demonstrate the application of our method using breast cancer patient data from Ontario, Canada. We validate this metric by showing it has a statistically significant association with survival. We also illustrate how it can be used to quantify population-level discordance in patient pathways relative to the recommended clinical pathways. We find that patients undertaking extra clinical activities constitute the major sources of discordance in most sections of the pathways. Further analysis is needed to understand the root of these discordances and identify opportunities for quality improvement.

\section*{Conflict of Interest}
The authors declare that they have no conflict of interest.

\newpage
\bibliographystyle{spbasic}      
\bibliography{references}

\newpage
\onecolumn
\appendix
\section*{Appendix}
\section{Pathway Requirements for Subgroups 3-9}
\label{app:tables}
\begin{ThreePartTable}

\begin{longtable}{llll}
\caption{Clinical pathway requirements for subgroup 3} \label{tab:activities3}\\

\toprule
Level 1          & Level 2                        & Level 3                      & Level 4                     \\
\toprule

Diagnosis (M)    & Diagnosis (M)    & Consultation (M)                  & Surgery Consultation (M)              \\
                 &                  &                              & Medical Oncologist Consultation (D)   \\
                 &                  &                              & Radiation Oncologist Consultation (D) \\
                 \cmidrule(l){3-4}
                 &                  & Screening Imaging (O)        & Screening mammogram (O)          \\
                 \cmidrule(l){3-4}
                 &                  & Diagnostic Imaging (M)       & Breast Ultrasound (O)            \\
                 &                  &                              & Ductogram (O)                    \\
\bottomrule
\end{longtable}
\begin{tablenotes}
        \scriptsize
        \item \textbf{M}: mandatory, \textbf{O}: optional, \textbf{D}: discordant, \textbf{AM}: alternative mandatory
    \end{tablenotes}
\end{ThreePartTable}

\begin{ThreePartTable}
\begin{longtable}{llll}
\caption{Clinical pathway requirements for subgroup 4} \label{tab:activities4}\\

\toprule
Level 1          & Level 2                        & Level 3                      & Level 4                     \\
\toprule
Diagnosis (M)    & Diagnosis (M)    & Consultation (M)                  & Surgery Consultation (M)              \\
                 &                  &                              & Medical Oncologist Consultation (O)   \\
                 &                  &                              & Radiation Oncologist Consultation (O) \\
                 \cmidrule(l){3-4}
                 &                  & Screening Imaging (O)        & Screening mammogram (O)          \\
                 \cmidrule(l){3-4}
                 &                  & Diagnostic Imaging (M)       & Breast Ultrasound (O)            \\
                 &                  &                              & Ductogram (O)                    \\
\bottomrule
\end{longtable}
\begin{tablenotes}
        \scriptsize
        \item \textbf{M}: mandatory, \textbf{O}: optional, \textbf{D}: discordant, \textbf{AM}: alternative mandatory, \textbf{AO}: alternative optional
    \end{tablenotes}
\end{ThreePartTable}

\begin{ThreePartTable}
\begin{longtable}{llll}
\caption{Clinical pathway requirements for subgroup 5} \label{tab:activities5}\\

\toprule
Level 1          & Level 2                        & Level 3                      & Level 4                     \\
\toprule
Diagnosis (M)    & Diagnosis (M)    & Consultation (M)                  & Surgery Consultation (M)              \\
                 &                  &                              & Medical Oncologist Consultation (D)   \\
                 &                  &                              & Radiation Oncologist Consultation (D) \\
                 \cmidrule(l){3-4}
                 &                  & Screening Imaging (O)        & Screening mammogram (O)          \\
                 \cmidrule(l){3-4}
                 &                  & Diagnostic Imaging (M)       & Breast Ultrasound (O)            \\
                 &                  &                              & Ductogram (O)                    \\
\bottomrule
\end{longtable}
\begin{tablenotes}
        \scriptsize
        \item \textbf{M}: mandatory, \textbf{O}: optional, \textbf{D}: discordant, \textbf{AM}: alternative mandatory
    \end{tablenotes}
\end{ThreePartTable}

\begin{ThreePartTable}
\begin{longtable}{llll}
\caption{Clinical pathway requirements for subgroup 6} \label{tab:activities6}\\

\toprule
Level 1          & Level 2                        & Level 3                      & Level 4                     \\
\toprule
Diagnosis (M)    & Diagnosis (M)    & Consultation (M)                  & Surgery Consultation (M)              \\
                 &                  &                              & Medical Oncologist Consultation (O)   \\
                 &                  &                              & Radiation Oncologist Consultation (O) \\
                 \cmidrule(l){3-4}
                 &                  & Screening Imaging (O)        & Screening mammogram (O)          \\
                 \cmidrule(l){3-4}
                 &                  & Diagnostic Imaging (M)       & Breast Ultrasound (O)            \\
                 &                  &                              & Ductogram (O)                    \\
\bottomrule
\end{longtable}
\begin{tablenotes}
        \scriptsize
        \item \textbf{M}: mandatory, \textbf{O}: optional, \textbf{D}: discordant, \textbf{AM}: alternative mandatory, \textbf{AO}: alternative optional
    \end{tablenotes}
\end{ThreePartTable}

\begin{ThreePartTable}
\begin{longtable}{llll}
\caption{Clinical pathway requirements for subgroup 7} \label{tab:activities7}\\

\toprule
Level 1          & Level 2                        & Level 3                      & Level 4                     \\
\toprule
Diagnosis (M)    & Diagnosis (M)    & Consultation (M)                  & Surgery Consultation (M)              \\
                 &                  &                              & Medical Oncologist Consultation (M)   \\
                 &                  &                              & Radiation Oncologist Consultation (D) \\
                 \cmidrule(l){3-4}
                 &                  & Screening Imaging (O)        & Screening mammogram (O)          \\
                 \cmidrule(l){3-4}
                 &                  & Diagnostic Imaging (M)       & Breast Ultrasound (O)            \\
                 &                  &                              & Ductogram (O)                    \\
\bottomrule
\end{longtable}
\begin{tablenotes}
        \scriptsize
        \item \textbf{M}: mandatory, \textbf{O}: optional, \textbf{D}: discordant, \textbf{AM}: alternative mandatory
    \end{tablenotes}
\end{ThreePartTable}

\begin{ThreePartTable}
\begin{longtable}{llll}
\caption{Clinical pathway requirements for subgroup 8} \label{tab:activities8}\\

\toprule
Level 1          & Level 2                        & Level 3                      & Level 4                     \\
\toprule
Diagnosis (M)    & Diagnosis (M)    & Consultation (M)                  & Surgery Consultation (M)              \\
                 &                  &                              & Medical Oncologist Consultation (M)   \\
                 &                  &                              & Radiation Oncologist Consultation (D) \\
                 \cmidrule(l){3-4}
                 &                  & Screening Imaging (O)        & Screening mammogram (O)          \\
                 \cmidrule(l){3-4}
                 &                  & Diagnostic Imaging (M)       & Breast Ultrasound (O)            \\
                 &                  &                              & Ductogram (O)                    \\
\bottomrule
\end{longtable}
\begin{tablenotes}
        \scriptsize
        \item \textbf{M}: mandatory, \textbf{O}: optional, \textbf{D}: discordant, \textbf{AM}: alternative mandatory
    \end{tablenotes}
\end{ThreePartTable}

\begin{ThreePartTable}
\begin{longtable}{llll}
\caption{Clinical pathway requirements for subgroup 9} \label{tab:activities9}\\

\toprule
Level 1          & Level 2                        & Level 3                      & Level 4                     \\
\toprule
Diagnosis (M)    & Diagnosis (M)    & Consultation (M)                  & Surgery Consultation (M)              \\
                 &                  &                              & Medical Oncologist Consultation (M)   \\
                 &                  &                              & Radiation Oncologist Consultation (D) \\
                 \cmidrule(l){3-4}
                 &                  & Screening Imaging (O)        & Screening mammogram (O)          \\
                 \cmidrule(l){3-4}
                 &                  & Diagnostic Imaging (M)       & Breast Ultrasound (O)            \\
                 &                  &                              & Ductogram (O)                    \\
\bottomrule
\end{longtable}
\begin{tablenotes}
        \scriptsize
        \item \textbf{M}: mandatory, \textbf{O}: optional, \textbf{D}: discordant, \textbf{AM}: alternative mandatory
    \end{tablenotes}
\end{ThreePartTable}

\section{Supplementary Figures and Tables}
\label{app:fig_table}

\begin{figure*}
    \includegraphics[width=0.6\textwidth]{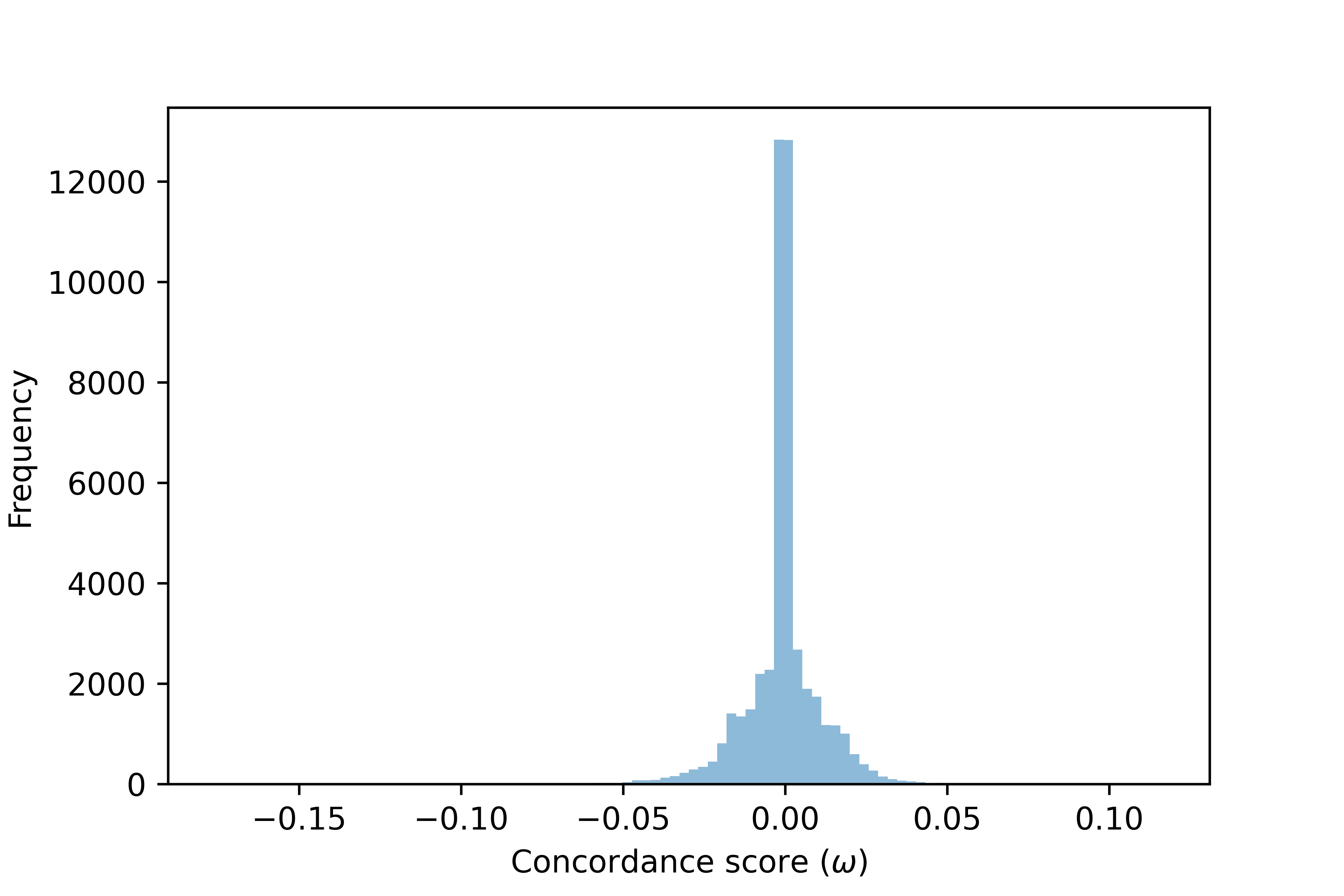}
    \caption{Distribution of difference in concordance scores between a bootstrapped sample and the original model}
    \label{fig:boot}
\end{figure*}

\begin{table*}
\caption{Risk of mortality associated with concordance score ($\omega$) for subgroups when patient data are swapped in model~\eqref{model:hierpatientref}}
\label{tab:omegasurvival-reversed}
\begin{tabular}{ccccccc}
\toprule
         & \multicolumn{3}{c}{Unadjusted}         & \multicolumn{3}{c}{Adjusted}           \\
        \cmidrule(r){2-4} \cmidrule(l){5-7}
Subgroup & HR   & 95\% CI      & p-value          & HR   & 95\% CI      & p-value          \\
\toprule
2        & 1.22 & (0.74, 2.02) & 0.4416 & 0.83 & (0.52, 1.33) & 0.4382 \\
3        & 0.77 & (0.40, 1.51) & 0.4524 & 0.60 & (0.31, 1.16) & 0.1291 \\
4        & 1.00 & (0.53, 1.88) & 1.0000 & 0.69 & (0.37, 1.29) & 0.2481 \\
5        & 1.12 & (0.66, 1.91) & 0.6681 & 1.05 & (0.61, 1.81) & 0.8556 \\
6        & 0.91 & (0.49, 1.70) & 0.7769 & 0.61 & (0.34, 1.11) & 0.1049 \\
7        & 1.53 & (0.63, 3.72) & 0.3433 & 1.34 & (0.57, 3.13) & 0.5033 \\
8        & 1.10 & (0.50, 2.42) & 0.8040 & 1.01 & (0.46, 2.21) & 0.9806 \\
9        & 1.09 & (0.55, 2.15) & 0.8103 & 0.72 & (0.37, 1.40) & 0.3387 \\
\bottomrule
\end{tabular}
\end{table*}

\end{document}